\documentclass[11pt]{article}
\usepackage{amsmath}
\usepackage{amsfonts}
\usepackage[mathscr]{euscript}
\usepackage{color}
\usepackage{hyperref}
\usepackage{soul}
\usepackage{graphics,graphicx}
\usepackage{subfig,float}
\newtheorem{theo}{Theorem}[section]
\newtheorem{lem}[theo]{Lemma}
\newtheorem{cor}[theo]{Corollary}

\newcommand{\mysection}[1]{\section{#1} \setcounter{equation}{0}}
\newcommand{\proof}{{\sc Proof.} \quad}
\newcommand{\proofc}{{\sc Proof} \ }
\newcommand{\be}{\begin{equation} \label}
\newcommand{\ee}{\end{equation}}
\newcommand{\bea}{\begin{eqnarray}\label}
\newcommand{\eea}{\end{eqnarray}}
\newcommand{\bas}{\begin{eqnarray*}}
\newcommand{\eas}{\end{eqnarray*}}
\newcommand{\bit}{\begin{itemize}}
\newcommand{\eit}{\end{itemize}}
\newcommand{\qed}{\hfill$\Box$ \vskip.2cm}
\newcommand{\nn}{\nonumber}
\newcommand{\norm}[1]{\|#1\|}
\newcommand{\R}{\mathbb{R}}
\newcommand{\N}{\mathbb{N}}
\newcommand{\pO}{\partial\Omega}
\newcommand{\bom}{\overline{\Omega}}
\newcommand{\eps}{\varepsilon}

\newcommand{\hra}{\hookrightarrow}
\newcommand{\io}{\int_\Omega}

\newcommand{\abs}{\\[5pt]}
\newcommand{\tm}{T_{max}}
\newcommand{\htm}{\widehat{T}_{max}}
\newcommand{\qp}{q^+(p)}
\newcommand{\qm}{q^-(p)}
\newcommand{\qpm}{q^\pm(p)}
\newcommand{\rs}{r_\star}
\newcommand{\ps}{p_\star}
\newcommand{\tp}{\widetilde{p}}
\begin{document}
\title{On the global existence and qualitative behavior of one-dimensional solutions to a model for urban crime}
\author{
Nancy Rodriguez\footnote{rodrign@colorado.edu}\\
{\small CU Boulder, Department of Applied Mathematics, Engineering Center, ECOT 225 }\\
{\small Boulder, CO 80309-0526}
\and
Michael Winkler\footnote{michael.winkler@math.uni-paderborn.de}\\
{\small Institut f\"ur Mathematik, Universit\"at Paderborn,}\\
{\small 33098 Paderborn, Germany} }

\date{}
\maketitle
\begin{abstract}
\noindent 
  We consider the no-flux initial-boundary value problem for the cross-diffusive evolution system
  \bas
   	\left\{ \begin{array}{ll}
	u_t = u_{xx} - \chi \big(\frac{u}{v} \partial_x u \big)_x - uv +B_1(x,t),
	\qquad & x\in \Omega, \ t>0, \\[1mm]
	v_t = v_{xx} +uv - v + B_2(x,t), 
	\qquad & x\in \Omega, \ t>0,
    	\end{array} \right.
  \eas
  which 
  was introduced by Short {\it et al.} in \cite{Short2008} with $\chi=2$ to describe the dynamics of urban crime.\abs
  In bounded intervals $\Omega\subset\R$ and with prescribed suitably regular nonnegative
  functions $B_1$ and $B_2$,
  we first prove the existence of global classical solutions for any choice of $\chi>0$ and all reasonably
  regular nonnegative initial data.  \abs
  We next address the issue of determining the qualitative behavior of solutions under appropriate assumptions
  on the asymptotic properties of $B_1$ and $B_2$.  
  Indeed, for arbitrary $\chi>0$ we obtain boundedness of the solutions given strict positivity of the average of $B_2$
  over the domain; 
  moreover, it is seen that imposing a mild decay assumption on $B_1$ implies that $u$ must
  decay to zero in the long-term limit.  
  Our final result, valid for all $\chi\in\left(0,\frac{\sqrt{6\sqrt{3}+9}}{2}\right),$ which contains
  the relevant value $\chi=2$, states
  that under the above decay assumption on $B_1$, if furthermore $B_2$ appropriately stabilizes to a
  nontrivial function $B_{2,\infty}$, then $(u,v)$ approaches the limit $(0,v_\infty)$, where $v_\infty$
  denotes the solution of 
  \bas
	\left\{ \begin{array}{l}
	-\partial_{xx}v_\infty + v_\infty = B_{2,\infty}, 
	\qquad x\in \Omega, \\[1mm]
	\partial_x v_{\infty}=0,
	\qquad x\in\pO.
	\end{array} \right.
  \eas
  We conclude with some numerical simulations exploring possible effects that may arise when 
  considering large values of $\chi$ not covered by our qualitative analysis.
  We observe that when $\chi$ increases,
  solutions may grow substantially on short time intervals, whereas only on large time scales
  diffusion will dominate and enforce equilibration.\abs
  {\bf Keywords:} urban crime, global existence, decay estimates, long-time behavior\\
  {\bf MSC (2010):} 35Q91 (primary); 35B40, 35K55 (secondary)
\end{abstract}
\newpage
\section{Introduction}\label{intro}
Driven by the need to understand the spatio-temporal dynamics of crime {\it hotspots}, which are regions in space
that have a disproportionately high level of crime, Short and collaborators introduced a reaction-advection-diffusion system to describe the evolution of urban crime in \cite{Short2008}.  
When posed in spatial one-dimensional domains $\Omega$, this system read
\be{0}
   	\left\{ \begin{array}{ll}
	u_t = u_{xx} - \chi \big(\frac{u}{v} v_x\big)_x - uv +B_1(x,t),
	\qquad & x\in \Omega, \ t>0, \\[1mm]
	v_t =  v_{xx} +uv - v + B_2(x,t), 
	\qquad & x\in \Omega, \ t>0,
    	\end{array} \right.
\ee
with the parameter $\chi$ fixed as
\be{chi}
	\chi=2,
\ee
and with given source functions $B_1$ and $B_2$.
In \eqref{0}, $u(x,t)$ represents the {\it density of criminal agents} and $v(x,t)$ the {\it attractiveness value}, which provides a measure
of how susceptible a certain location $x$ is to crime at time $t.$
System \eqref{0} was derived from an agent-based model rooted on the assumption of ``routine activity theory",  
a criminology theory stating that opportunity is the most important factor leading to crime  \cite{Cohen1979, Felson1987}.   The system models two sociological effects: the `repeat and
near-repeat victimization' effect and the `broken-windows theory'.  The former has been observed in residential burglary data and alludes to the 
fact that the burglarization of a house increases the probability of that same house, as well as neighboring houses, to be burgled again within a short period
of time following the original burglary \cite{Johnson1997, Short2009}.  The latter is the theory that, in a sense, crime is self-exciting - crime tends to lead to more crime \cite{Kelling1982}.\abs
From the first equation in \eqref{0} we see that criminal agents move according to a combination of conditional and unconditional diffusion.  
The conditional diffusion is a biased movement toward high concentrations of the attractiveness value, which leads to the taxis term seen
in the first equation.  We stress that the coefficient $\chi = 2$ in front of the taxis term, which we shall see adds a challenge, comes from the first principles derivation of system \eqref{0}
and thus it is important that our theory cover this case -- see \cite{Short2008} for more details.
The assumption that criminal agents 
abstain from committing a second crime leads to decay term $-uv.$  
Indeed, roughly 
speaking, the expected number of crime is given by $uv$ and so the expected number of criminal agents removed is $uv$.  
The prescribed non-negative term $B_1(x)$ describes the 
introduction of criminal agents into the system.  Furthermore, the repeat victimization effect assumes that each criminal activity increases the attractiveness value leading to the $+uv$ term in the second equation of 
\eqref{0}, while the near-repeat victimization effect leads to the unconditional diffusion also observed in that equation.  
Finally, the assumption that certain neighborhoods tend to be more crime-prone than others, whatever these reasons may be,
is included in the prescribed non-negative term $B_2(x).$\abs
The introduction of system \eqref{0} has generated a great deal of activity related to the analysis of \eqref{0}, which have contributed to the mathematical theory
as well as to the understanding of crime dynamics.  
For example, the emergence and suppression of {\it hotspots} was studied by Short {\it et al.} in \cite{Short2010}, 
providing insight into the effectiveness of {\it hotspot policing}.    
The existence and stability of localized patterns representing hotspots has been studied in 
various works -- see  
 \cite{Berestycki2014a, Cantrell2012, Gu2016, Kolokolnikiv2014,Tse2008}. 
A more general class of systems was proposed for the dynamics of criminal activity by Berestycki and Nadal in \cite{Berestycki2010} --
see also \cite{Berestycki2013c} for an analysis of these models.  The system \eqref{0} has also been generalized in various directions.  For example,
the incorporation of law enforcement has been proposed and analyzed in \cite{Jones2010, Ricketson2010, Zipkin2014}; the movement of 
{\em commuter} criminal agents was modeled in \cite{Chaturapruek2013} through the use of L\'evy flights.  The dynamics of crime 
has also be studied with the use of dynamics systems, we refer the readers to \cite{McMillon2014, Nuno2011}.  It is also important to 
note that the work in \cite{Short2008} has been the impetus for the use of PDE type models to gain insight into various other 
social phenomena -- see for example \cite{Barbaro2013, Rodriguez2014a, Smith2012}.  Interested readers are referred to 
the comprehensive review of mathematical models and theory for criminal activity in \cite{DOrsogna2015}.\abs
From a perspective of mathematical analysis, (\ref{0}) shares essential ingredients with the celebrated
Keller-Segel model for chemotaxis processes in biology, which in its simplest form can be obtained on 
considering the constant sensitivity function $S\equiv 1$ in
\be{KS}
   	\left\{ \begin{array}{ll}
	u_t = \Delta u - \nabla \cdot (uS(v)\nabla v),
	\qquad & x\in \Omega, \ t>0, \\[1mm]
	v_t =  \Delta v - v + u,
	\qquad & x\in \Omega, \ t>0.
    	\end{array} \right.
\ee
Here the interplay of such cross-diffusive terms
with the linear production mechanism expressed in the second equation is known to have a strongly destabilizing
potential in multi-dimensional situations: when posed under no-flux boundary conditions 
in bounded domains $\Omega\subset\R^n$, $n\ge 1$, (\ref{KS}) is globally well-posed in the case $n=1$
(\cite{osaki_yagi}), whereas some solutions may blow up in finite time
when either $n= 2$ and the conserved quantity $\io u(\cdot,t)$ is suitably large 
(\cite{herrero_velazquez}), or when $n\ge 3$ (\cite{win_JMPA}; cf.~also the recent survey
\cite{BBTW}).\abs
That, in contrast to this, decaying sensitivities may exert a substantial regularizing effect is indicated by
the fact that if 
e.g.~$S(v)=\frac{a}{(1+bv)^\alpha}$ for all $v\ge 0$ and some $a>0$, $b>0$ and $\alpha>1$, then actually for
arbitrary $n\ge 1$ global bounded solutions to (\ref{KS}) always exist (\cite{win_MANA}).
However, in the particular case of the so-called logarithmic sensitivity given by $S(v)=\frac{\chi}{v}$ for $v>0$ with
$\chi>0$, as present in (\ref{0}), the situation seems less clear in that global bounded solutions so far have been
constructed only under smallness conditions of the form $\chi<\sqrt{\frac{2}{n}}$ (\cite{biler}, 
\cite{win_chemosing}), with a slight extension
up to the weaker condition $\chi<\chi_0$ with some $\chi_0\in (1.015,2)$ possible when $n=2$ (\cite{lankeit_MMAS}); 
for larger values of $\chi$ including the choice in (\ref{chi}), in the case $n\ge 2$
only certain global weak solutions to (\ref{KS}),
possibly becoming unbounded in finite time, are known to exist in various generalized frameworks
(\cite{win_chemosing}, \cite{stinner_win}, \cite{lw2}).\abs
With regard to issues of regularity and boundedness, the situation in (\ref{0}) seems yet more delicate than
in the latter version of (\ref{KS}): In (\ref{0}), namely, the production of the attractiveness value occurs
in a nonlinear manner, which in comparison to (\ref{KS}) may further stimulate the self-enhanced generation
of large cross-diffusive gradients.
To the best of our knowledge, no results on global existence have
been found so far for any version of (\ref{KS}) in which such reaction terms is introduced, even in spatially
one-dimensional cases, and it seems far from obvious to which extent such mechanisms can be compensated by
the supplementary absorptive term $-uv$ in the first equation of (\ref{0}).\abs
Accordingly, the literature on initial-value problems for (\ref{0}) is still at quite an early stage
and actually limited to a first local existence and uniqueness result achieved in \cite{Rodriguez2010}.
Statements on global existence have been obtained only for certain modified versions which contain additional
regularizing ingredients (\cite{Manasevich2012}, \cite{Rodriguez2013}). \abs
{\bf Main results.} \quad
In the present work we attempt to undertake a first step into a qualitative theory for the full original model
from \cite{Short2008} by developing an approach capable of analyzing 
the spatially one-dimensional system (\ref{0}) in a range of parameters including the choice given in (\ref{chi}).
Here we will first concentrate be on establishing a result on global existence of classical solutions under mild 
assumptions on $\chi$, $B_1$ and $B_2$. Our second focus will be on the derivation of qualitative
solution properties under additional assumptions.\abs
In order to specify the setup for our analysis, for a given parameter $\chi>0$
let us consider (\ref{0}) along with the boundary conditions
\be{0b}
	u_x=v_x=0,
	\qquad x\in \pO, \ t>0,
\ee
and the initial conditions
\be{0i}
	u(x,0)=u_0(x), \quad v(x,0)=v_0(x), 
	\qquad  x\in\Omega,
\ee
in a bounded open interval $\Omega\subset \R$. We assume throughout the sequel that
\be{B}
	\mbox{$B_1$ and $B_2$ are nonnegative bounded functions belonging to 
	$C^\vartheta_{loc}(\bom\times [0,\infty))$ for some $\vartheta \in (0,1)$,}
\ee
and that 
\be{init}
	\left\{ \begin{array}{ll}
	u_0\in C^0(\bom),&
	\quad \mbox{with } u_0\ge 0 \mbox{ in } \Omega,  \\[1mm]
	v_0\in W^{1,\infty}(\Omega),&
	\quad \mbox{with } v_0> 0 \mbox{ in } \bom.
	\end{array} \right.
\ee
In this general framework, we shall see that in fact for arbitrary $\chi>0$, the problem 
\eqref{0}, \eqref{0b}, \eqref{0i} 
is globally well-posed in the following sense.
\begin{theo}\label{theo5}
  Let $\chi>0$ and suppose that $B_1$ and $B_2$ satisfy \eqref{B}.
  Then for any choice of $u_0$ and $v_0$ fulfilling \eqref{init}, the problem 
  \eqref{0}, \eqref{0b}, \eqref{0i} 
  possesses a global classical
  solution, for each $r>1$ uniquely determined by the inclusions
  \be{5.1}
	\left\{ \begin{array}{l}
	u\in C^0(\bom\times [0,\infty)) \cap C^{2,1}(\bom\times (0,\infty)), \\[1mm]
	v\in C^0([0,\infty);W^{1,r}(\Omega)) \cap C^{2,1}(\bom\times (0,\infty)),
	\end{array} \right.
  \ee
  for which $u,v>0$ in $\bom\times (0,\infty)$.
\end{theo}
The qualitative behavior of these solutions, especially on large time scales, will evidently depend on 
respective asymptotic properties of the parameter functions $B_1$ and $B_2$.
Our efforts in this direction will particularly make use of either suitable assumptions on large-time 
decay of $B_1$ or of certain weak but temporally uniform positivity properties of $B_2$.
Specifically, in our analysis we will alternately refer to the hypotheses
\be{H1}
	\int_0^\infty \io B_1 < \infty,
	\tag{H1}
\ee
and, in a weaker form,
\be{H1'}
	\int_t^{t+1} \io B_1(x,s) dxds \to 0
	\qquad \mbox{as } t\to\infty,
	\tag{H1'}
\ee
on decay of $B_1$, and
\be{H2}
	\inf_{t>0} \io B_2(x,t)dx > 0,
	\tag{H2}
\ee
on the positivity of $B_2$.
In some places we will also assume that $B_2$ stabilizes in the sense that
\be{H3}
	\int_t^{t+1} \io \Big(B_2(x,s)-B_{2,\infty}(x)\Big)^2 dxds \to 0
	\qquad \mbox{as } t\to\infty,
	\tag{H3}
\ee
holds with some $B_{2,\infty}\in L^2(\Omega)$.\abs
Indeed, the assumption \eqref{H2} implies boundedness of both solution components, and under the additional
requirement that \eqref{H1'} be valid, $u$ must even decay in the large time limit.
\begin{theo}\label{theo51}
  Let $\chi>0$ and suppose that \eqref{B} and \eqref{init} are fulfilled.
  If moreover \eqref{H2} holds, then there exists $C>0$ with the property that the solution $(u,v)$ of
  \eqref{0}, \eqref{0b}, \eqref{0i} 
  satisfies
  \be{51.1}
	u(x,t) \le C
	\qquad \mbox{for all $x\in\Omega$ and } t>0,
  \ee
  and
  \be{51.2}
	\frac{1}{C} \le v(x,t) \le C
	\qquad \mbox{for all $x\in\Omega$ and } t>0.
  \ee
  If additionally \eqref{H1'} is valid, then
  \be{51.3}
	u(\cdot,t) \to 0
	\quad \mbox{in } L^\infty(\Omega)
	\qquad \mbox{as } t\to\infty.
  \ee
\end{theo}
We shall secondly see that for all $\chi$ within an appropriate range, including the relevant value $\chi=2$,
also the mere assumption \eqref{H1} is sufficient for boundedness, at least of the second solution component,
and that moreover the latter even stabilizes when additionally \eqref{H3} is satisfied.
\begin{theo}\label{theo52}
  Let $\chi>0$ be such that
  \be{52.1}
	\chi < \frac{\sqrt{6\sqrt{3}+9}}{2}=2.201834...,
  \ee
  and let $B_1$ and $B_2$ be such that besides \eqref{B}, also \eqref{H1} holds.
  Then for each pair $(u_0,v_0)$ fulfilling \eqref{init}, one can find $C>0$ such that
  the solution $(u,v)$ of 
  \eqref{0}, \eqref{0b}, \eqref{0i} 
  satisfies
  \be{52.2}
	v(x,t) \le C,
	\qquad \mbox{for all $x\in\Omega$ and } t>0.
  \ee
  Furthermore, if \eqref{H3} is valid with some $B_{2,\infty}\in L^2(\Omega)$, then 
  \be{52.3}
	v(\cdot,t) \to v_\infty
	\quad \mbox{in } L^\infty(\Omega)
	\qquad \mbox{as } t\to\infty,
  \ee
  where $v_\infty$ denotes the solution to the boundary value problem
  \be{vinfty}
	\left\{ \begin{array}{l}
	-\partial_{xx}v_\infty + v_\infty = B_{2,\infty}, 
	\qquad x\in \Omega, \\[1mm]
	\partial_xv_{\infty}=0,
	\qquad x\in\pO.
	\end{array} \right.
  \ee
\end{theo}
Let us finally state an essentially immediate consequence of Theorem \ref{theo51} and Theorem \ref{theo52} 
under slightly sharper but yet quite practicable assumptions.
\begin{cor}\label{cor53}
  Let $\chi\in (0,\frac{\sqrt{6\sqrt{3}+9}}{2})$, and suppose that the functions $B_1$ and $B_2$ are such that
  beyond \eqref{B} and \eqref{H1} we have
  \be{53.2}
	B_2(\cdot,t) \to B_{2,\infty}
	\quad \mbox{a.e.~in } \Omega
	\qquad \mbox{as $t\to\infty$},
  \ee
  with some $0\not\equiv B_{2,\infty} \in L^1(\Omega)$. Then for each $u_0$ and $v_0$ satisfying \eqref{init},
  the corresponding solution $(u,v)$ of
  \eqref{0}, \eqref{0b}, \eqref{0i} 
  has the properties that
  \bas
	u(\cdot,t) \to 0
	\quad \mbox{in } L^\infty(\Omega)
	\qquad \mbox{as } t\to\infty,
  \eas
  and
  \bas
	v(\cdot,t) \to v_\infty
	\quad \mbox{in } L^\infty(\Omega)
	\qquad \mbox{as } t\to\infty,
  \eas
  where $v_\infty$ solves \eqref{vinfty}.
\end{cor}
\proof		
  In view of the dominated convergence theorem, \eqref{53.2} along with the boundedness of $B_2$ entails that
  actually $B_{2,\infty}\in L^\infty(\Omega)$, that \eqref{H3} holds and that moreover
  $\io B_2(\cdot,t)\to \io B_{2,\infty}\ne 0$ as $t\to\infty$, whence for some $t_0>0$ we have
  $\inf_{t>t_0} \io B_2(\cdot,t)>0$.
  The claim therefore results on applying Theorem \ref{theo51} and Theorem \ref{theo52} with $(u,v,B_1,B_2)(x,t)$
  replaced by $(u,v,B_1,B_2)(x,t_0+t)$ for $(x,t)\in\bom\times [0,\infty)$.
\qed
{\bf Outline.} \quad
After asserting local existence of solutions and some of their basic features in
Section \ref{sec:local},
in Section \ref{sec:fund} we will derive some fundamental estimates resulting from an analysis of the coupled
functional $\io u^p v^q$ which indeed enjoys a certain entropy-type property if, in dependence on the size of $\chi$,
the crucial exponent $p$ therein is small enough and $q$ belongs to an appropriate range.
Accordingly implied consequences on regularity features will thereafter enable us to verify 
Theorem \ref{theo5} and Theorem \ref{theo51} in Section \ref{sec:global}. 
Finally, Section \ref{sec:longtime} will contain our proof of Theorem \ref{theo52}, where we highlight already here that
particular challenges will be linked
to the derivation of $L^\infty$ bounds for $v$, and that these will be accomplished on the basis of a recursive argument
available under the assumption (\ref{52.1}).

\mysection{Local existence and basic estimates}\label{sec:local}
Let us first make sure that our overall assumptions warrant local-in-time solvability of
\eqref{0}, \eqref{0b}, \eqref{0i}, along with a convenient extensibility criterion.
\begin{lem}\label{lem_loc}
  Under the assumptions of Theorem \ref{theo5}, there exist $\tm\in (0,\infty]$ and a uniquely determined pair $(u,v)$
  of functions
  \bas
	\left\{ \begin{array}{l}
	u \in C^0(\bom\times [0,\tm)) \cap C^{2,1}(\bom\times (0,\tm)), \\[1mm]
	v \in \bigcap\limits_{r>1} C^0([0,\tm); W^{1,r}(\Omega)) \cap C^{2,1}(\bom\times (0,\tm)),
	\end{array} \right.
  \eas
  which solve 
  \eqref{0}, \eqref{0b}, \eqref{0i} 
  classically in $\bom\times [0,\tm).$ Moreover, $u>0$ and $v>0$ in $\bom\times (0,\tm)$ and 
  \be{ext_crit}
	\mbox{either $\tm=\infty$, \quad or \quad}
	\limsup_{t\nearrow\tm} \Big\{ \|u(\cdot,t)\|_{L^\infty(\Omega)}
	+ \Big\|\frac{1}{v(\cdot,t)}\Big\|_{L^\infty(\Omega)}
	+ \|v_x(\cdot,t)\|_{L^r(\Omega)} \Big\} = \infty
	\quad \mbox{for all } r>1.
  \ee
\end{lem}
\proof
  The results is a straightforward application of well-established techniques from the theory of tridiagonal
  cross-diffusive systems (\cite{amann}, specifically applied to chemotaxis systems \cite{horstmann_win}).
\qed
Throughout the sequel, without explicit further mentioning we shall assume the requirements of Theorem
\ref{theo5} to be met, and let $u, v$ and $\tm$ be as provided by Lemma \ref{lem_loc}.\abs
In order to derive some basic features of this solution, let us recall 
the following well-known 
pointwise positivity property of the Neumann heat semigroup $(e^{t\Delta})_{t\ge 0}$ on the bounded
real interval $\Omega$ (cf.~e.g.~\cite[Lemma 3.1]{hpw}).
\begin{lem}\label{lem24}
  Let $\tau>0$. Then there exists a constant $C>0$ such that for all nonnegative $\varphi\in C^0(\bom)$,
  \bas
	e^{t\Delta} \varphi \ge C \io \varphi
	\quad \mbox{in } \Omega
	\qquad \mbox{for all } t> \tau.
  \eas
\end{lem}
Using the previous lemma along with a parabolic comparison argument, we obtain a basic but important pointwise lower estimate
for the second solution component.  This lower bound is local-in-time for arbitrary $B_1$ and $B_2$ and global-in-time when
\eqref{H2} is satisfied.
\begin{lem}\label{lem001}
  For all $T>0$ there exists $C(T)>0$ such that with $\tm$ from Lemma \ref{lem_loc}, for $\htm:=\min\{T,\tm\}$ we have
  \be{001.1}
	v(x,t) \ge C(T),
	\qquad \mbox{for all $x\in\Omega$ and } t\in (0,\htm),
  \ee
  with  \be{001.2}
	\inf_{T>0} C(T)>0,
	\qquad \mbox{if \eqref{H2} is valid.}
  \ee
\end{lem}
\proof
  We represent $v$ according to
  \bea{001.3}
	v(\cdot,t)= e^{t(\Delta-1)} v_0
	+ \int_0^t e^{(t-s)(\Delta-1)} u(\cdot,s) v(\cdot,s) ds
	+ \int_0^t e^{(t-s)(\Delta-1)} B_2(\cdot,s) ds,
	\qquad t\in (0,\tm),
  \eea
  and observe that here by the comparison principle for the Neumann problem associated with the heat equation,
  the second summand on the right is nonnegative, whereas 
  \be{001.33}
	e^{t(\Delta-1)} v_0
	\ge \Big\{ \inf_{x\in\Omega} v_0(x) \Big\} \cdot e^{-t}
	\qquad \mbox{for all } t>0.
  \ee
  To gain a pointwise lower estimate for the rightmost integral in \eqref{001.3}, we invoke Lemma \ref{lem24}
  to find $c_1>0$ such that with $\tau:=\min\{1,\frac{1}{3}\tm\}$, for any nonnegative $\varphi\in C^0(\bom)$ we have
  \bas
	e^{t\Delta}\varphi \ge c_1 \io \varphi
	\quad \mbox{in } \Omega
	\qquad \mbox{for all } t>\frac{\tau}{2},
  \eas
  which implies that
  \bas
	\int_0^t e^{(t-s)(\Delta-1)} B_2(\cdot,s) ds
	&\ge& \int_0^{t-\frac{\tau}{2}} e^{-(t-s)} e^{(t-s)\Delta} B_2(\cdot,s) ds \\
	&\ge& \int_0^{t-\frac{\tau}{2}} e^{-(t-s)} \cdot \Big\{ c_1 \io B_2(\cdot,s) \Big\} ds \\
	&\ge& c_1 c_2 \int_0^{t-\frac{\tau}{2}} e^{-(t-s)} ds \\
	&=& c_1 c_2 \cdot \Big(e^{-\frac{\tau}{2}} - e^{-t}\Big) \\
	&\ge& c_3:= c_1 c_2 \cdot \Big( e^{-\frac{\tau}{2}} - e^{-\tau}\Big)
	\qquad \mbox{for all } t>\tau
  \eas
  with $c_2:=\inf_{t>0} \io B_2(\cdot,t) \ge 0$.
  Together with \eqref{001.33} and \eqref{001.3}, this entails that
  \bas
	v(\cdot,t) \ge \Big\{ \inf_{x\in\Omega} v_0(x) \Big\} \cdot e^{-T} + c_3
	\quad \mbox{in } \Omega
	\qquad \mbox{for all } t\in (\tau,\htm),
  \eas
  and that
  \bas
	v(\cdot,t) \ge \Big\{ \inf_{x\in\Omega} v_0(x) \Big\} \cdot e^{-\tau}
	\quad \mbox{in } \Omega
	\qquad \mbox{for all } t\in (0,\tau],
  \eas
  and thereby establishes both \eqref{001.1} and \eqref{001.2}.
\qed
Further fundamental properties of \eqref{0} are connected to the evolution of the total mass $\io u$
and the associated total absorption rate $\io uv$.
We formulate these properties in such a way that important dependences of the appearing constants are accounted for
in order to provide statements that will be useful for our asymptotic
analysis in Theorem \ref{theo51} and Theorem \ref{theo52}.  %
\begin{lem}\label{lem01}
  For all $T>0$ there exists $C(T)>0$ such that with $\htm:=\min\{T,\tm\}$,	
  \be{01.1}
	\io u(\cdot,t) \le C(T)
	\qquad \mbox{for all } t\in (0,\htm),
  \ee
  where
  \be{01.11}
	\sup_{T>0} C(T) < \infty
	\qquad \mbox{if either \eqref{H1} or \eqref{H2} hold.}
  \ee
  Moreover, for all $T>0$ and each $\xi\in (0,T)$ there exists $K(T,\xi)>0$ with the properties that
  \be{01.2}
	\int_t^{t+\xi} \io uv \le K(T,\xi)
	\qquad \mbox{for all } t\in (0,\htm-\xi)
  \ee
  and
  \be{01.21}
	\sup_{T>\xi} K(T,\xi) < \infty 
	\quad \mbox{for all } \xi>0
	\qquad \mbox{if \eqref{H2} holds}
  \ee
  as well as
  \be{01.22}
	\sup_{T>0} \sup_{\xi\in (0,T)} K(T,\xi) < \infty
	\qquad \mbox{if \eqref{H1} holds.}
  \ee
\end{lem}
\proof
  Integrating the first equation in (\ref{0}) yields
  \be{01.3}
	\frac{d}{dt} \io u = - \io uv + \io B_1
	\qquad \mbox{for all } t\in (0,\tm)
  \ee
  and hence
  \be{01.4}
	\io u(\cdot,t) \le c_1(T):=\io u_0 + \int_0^T \io B_1
	\qquad \mbox{for all } t\in (0,\htm)
  \ee
  as well as
  \be{01.5}
	\int_t^{t+\xi} \io uv \le \io u(\cdot,t) + \int_t^{t+\xi} \io B_1 \le c_1(T) + c_1(2T)
	\qquad \mbox{for all $t \in (0,\htm-\xi)$ and any } \xi \in (0,T).
  \ee
  For general $B_1$ and $B_2$, (\ref{01.4}) and (\ref{01.5}) directly imply (\ref{01.1}) and (\ref{01.2}) with
  $C(T):=c_1(T)$ and $K(T,\xi):=c_1(T)+c_1(2T)$ for $T>0$ and $\xi\in (0,T)$, 
  and if in addition (\ref{H1}) holds, then $c_1(T)\le c_2:=\int_0^\infty \io B_1$ for all $T>0$ and thus 
  (\ref{01.4}) and (\ref{01.5}) moreover show that $C(T)\le c_2$ in this case.\abs
  Assuming the hypothesis (\ref{H2}) henceforth, we recall that thanks to the latter, Lemma \ref{lem001} implies the
  existence of $c_3>0$ fulfilling
  $v\ge c_3$ in $\Omega\times (0,\tm)$, whence going back to (\ref{01.3}) we see that then
  \be{01.6}
	\frac{d}{dt} \io u
	+ \frac{1}{2} \io uv \le - \frac{c_3}{2} \io u + c_4
	\qquad \mbox{for all } t\in (0,\tm)
  \ee
  with $c_4:=|\Omega|\cdot\|B_1\|_{L^\infty(\Omega\times (0,\infty))}$. 
  By an ODE comparison, this firstly ensures that
  \bas
	\io u(\cdot,t) \le c_5:=\max \Big\{ \io u_0, \frac{2c_4}{c_3}\Big\}
	\qquad \mbox{for all } t\in (0,\tm),
  \eas
  whereupon an integration in (\ref{01.6}) shows that furthermore
  \bas
	\frac{1}{2} \int_t^{t+\xi} \io uv 
	\le \io u(\cdot,t) + c_4\xi
	\le c_5 + c_4
	\qquad \mbox{for all } t\in (0,\tm-\xi)
  \eas
  and that hence indeed the estimates in (\ref{01.1}) and (\ref{01.2}) can actually be achieved to be independent of $T$
  also when (\ref{H2}) holds.
\qed
The previous lemma has the following consequence for the time evolution of $\io v$.
\begin{lem}\label{lem02}
  For all $T>0$ there exists $C(T)>0$ such that with $\htm:=\min\{T,\tm\}$ and $\tau:=\min\{1,\frac{1}{3}\tm\}$ we have
  \be{02.1}
	\io v(\cdot,t) \le C(T)
	\qquad \mbox{for all } t\in (0,\htm)
  \ee
  and 
  \be{02.2}
	\sup_{T>0} C(T) < \infty
	\qquad \mbox{if either \eqref{H1}	or \eqref{H2} holds.}
  \ee
\end{lem}
\proof
  From the second equation in (\ref{0}) we obtain that
  \be{02.3}
	\frac{d}{dt} \io v + \io v = \io uv + \io B_2
	\qquad \mbox{for all } t\in (0,\tm).
  \ee
  Here we only need to observe that thanks to 
  Lemma \ref{lem01} and the boundedness of $B_2$ we can find
  $c_1(T)>0$ such that for $h(t):=\io u(\cdot,t)v(\cdot,t) + \io B_2(\cdot,t)$, 
  $t\in (0,\tm)$, we have
  \bas
	\int_t^{t+\tau} h(s) ds \le c_1(T)
	\qquad \mbox{for all } t\in (0,\htm-\tau),
  \eas
  and that
  \bas
	\sup_{T>0} c_1(T)<\infty
	\qquad \mbox{if either (\ref{H1}) or (\ref{H2}) hold.}
  \eas
  Therefore, extending $h$ by zero to all of $(0,\infty)$ we may apply Lemma \ref{lem_ssw} from the appendix below
  so as to derive (\ref{02.1}) and (\ref{02.2}) from (\ref{02.3}).
\qed
\mysection{Fundamental estimates resulting from an analysis of $\io u^p v^q$} \label{sec:fund}
The main goal of this section consists of deriving spatio-temporal $L^2$ bounds for both
$u_x$ and $v_x$ with appropriate solution-dependent weight functions.
This will be accomplished in Lemma \ref{lem1} through an analysis of the functional $\io u^p v^q$
for adequately small $p\in (0,1)$ and certain positive $q<1-p$ taken from a suitable interval.
Entropy-like properties of 
functionals containing multiplicative couplings of both solution components have played important
roles in the analysis of several chemotaxis problems at various stages of existence and regularity
theory, but in most precedent cases the respective dependence on the unknown is 
either of strictly convex type with respect to both solution components separately 
(\cite{tao_small}, \cite{taowin_consumption}, \cite{win_ARMA}, \cite{win_MANA}), 
or at least exhibits some superlinear growth with respect to the full solution couple 
when viewed as a whole (\cite{biler}, \cite{Manasevich2012}).
In addition, contrary to related situations addressing singular sensitivities of the form in (\ref{0})
(\cite{win_chemosing}, \cite{stinner_win}), the additional zero-order nonlinearities $uv$ appearing in the
present context of (\ref{0}) will
require adequately coping with respectively occurring superlinear terms (cf.~e.g.~(\ref{1.11}) below).
In preparation to a corresponding testing procedure, we will therefore independently derive a regularity
property of $v$ by using a quasi-entropy property of the functional $-\io v^q$ for arbitrary $q\in (0,1)$.
\subsection{A spatio-temporal bound for $v$ in $L^r$ for $r<3$}
By means of a standard testing procedure solely involving the second equation in (\ref{0}),
thanks to Lemma \ref{lem02} and the nonnegativity of $B_2$ we can derive the following.
\begin{lem}\label{lem21}
  Let $q\in (0,1)$. Then for each $T>0$ one can find $C(T)>0$ with the properties that
  \be{21.1}
	\int_t^{t+\tau} \io v^{q-2} v_x^2 \le C(T)
	\qquad \mbox{for all } t\in (0,\htm-\tau)
  \ee
  and that
  \be{21.11}
	\sup_{T>0} C(T)<\infty
	\qquad \mbox{if either \eqref{H1} or \eqref{H2} hold,}
  \ee
  where again $\htm:=\min\{T,\tm\}$ and $\tau:=\min\{1,\frac{1}{3}\tm\}$.
\end{lem}
\proof
  As $v>0$ in $\bom\times [0,\tm)$ by Lemma \ref{lem001}, we may test the second equation in (\ref{0}) by $v^{q-1}$
  to see that
  \bas
	\frac{1}{q} \frac{d}{dt} \io v^q
	&=& (1-q) \io v^{q-2} v_x^2
	+ \io uv^q - \io v^q + \io B_2 v^{q-1} \\
	&\ge& (1-q)\io v^{q-2} v_x^2 - \io v^q
	\qquad \mbox{for all } t\in (0,\tm),
  \eas
  which on further integration yields that
  \be{21.2}
	(1-q)\int_t^{t+\tau} \io v^{q-2} v_x^2
	\le \frac{1}{q} \io v^q(\cdot,t+\tau) + \int_t^{t+\tau} \io v^q
	\qquad \mbox{for all } t\in (0,\tm-\tau).
  \ee
  Since with $c_1:=|\Omega|^{1-q}$ we have
  \bas
	\io v^q \le c_1 \bigg\{ \io v\bigg\}^q
	\qquad \mbox{for all } t\in (0,\tm)
  \eas
  by the H\"older inequality, from (\ref{21.2}) we this obtain that
  \bas
	(1-q)\int_t^{t+\tau} \io v^{q-2} v_x^2
	\le \Big(\frac{1}{q}+1\Big) \cdot c_1 \cdot \bigg\{ \sup_{s\in (0,\htm)} \io v(\cdot,s) \bigg\}^q
	\qquad \mbox{for all } t\in (0,\htm-\tau),
  \eas
  which in view of Lemma \ref{lem02} implies (\ref{21.1}) and (\ref{21.11}).
\qed
Thanks to the fact that the considered spatial setting is one-dimensional, an interpolation of the above
result with the outcome of Lemma \ref{lem02} has a natural consequence on space-time integrability of $v$.
\begin{lem}\label{lem22}
  Given $r\in (1,3)$, for any $T>0$ one can fix $C(T)>0$ such that
  \be{22.1}
	\int_t^{t+\tau} \io v^r \le C(T)
	\qquad \mbox{for all } t\in (0,\htm-\tau)
  \ee
  and that
  \be{22.11}
	\sup_{T>0} C(T)<\infty
	\qquad \mbox{if either \eqref{H1} or \eqref{H2} holds,}
  \ee
  where $\htm:=\min\{T,\tm\}$ and $\tau:=\min\{1,\frac{1}{3}\tm\}$.
\end{lem}
\proof
  We may assume that $r\in (2,3)$ and then let $q:=r-2\in (0,1)$ to obtain from Lemma \ref{lem21} that there exists
  $c_1(T)>0$ such that
  \be{22.2}
	\int_t^{t+\tau} \io [(v^\frac{q}{2})_x]^2 \le c_1(T)
	\qquad \mbox{for all } t\in (0,\htm-\tau),
  \ee
  while Lemma \ref{lem02} provides $c_2(T)>0$ fulfilling
  \be{22.3}
	\|v^\frac{q}{2}\|_{L^\frac{2}{q}(\Omega)}^\frac{2}{q} = \io v \le c_2(T)
	\qquad \mbox{for all } t\in (0,\htm),
  \ee
  where
  \be{22.33}
	\sup_{T>0} \Big(c_1(T)+ c_2(T)\Big) < \infty
	\qquad \mbox{if either (\ref{H1}) or (\ref{H2}) holds.}
  \ee
  Now, from the Gagliardo-Nirenberg inequality we know that there exists $c_3>0$ satisfying
  \bas
	\int_t^{t+\tau} \io v^r
	&=& \int_t^{t+\tau} \|v^\frac{q}{2}(\cdot,s)\|_{L^\frac{2r}{q}(\Omega)}^\frac{2r}{q} ds \\
	&\le& c_3 \int_t^{t+\tau} \bigg\{
	\Big\| (v^\frac{q}{2})_x(\cdot,s) \Big\|_{L^2(\Omega)}^\frac{2(r-1)}{q+1}
	\|v^\frac{q}{2}(\cdot,s)\|_{L^\frac{2}{q}(\Omega)}^\frac{2(q+r)}{q(q+1)}
	+ \|v^\frac{q}{2}(\cdot,s)\|_{L^\frac{2}{q}(\Omega)}^\frac{2r}{q} \bigg\} ds
  \eas
  for all $t\in (0,\htm-\tau)$, so that since 
  \bas
	\frac{2(r-1)}{q+1}=2
	\qquad \mbox{and} \qquad
	\frac{2(q+r)}{q(q+1)} = \frac{4}{r-2}=\frac{4}{q},
  \eas
  due to our choice of $q$ we obtain from (\ref{22.2}) and (\ref{22.3}) that
  \bas
	\int_t^{t+\tau} \io v^r
	&\le& c_3 \int_t^{t+\tau} \bigg\{ \Big\|(v^\frac{q}{2})_x(\cdot,s)\Big\|_{L^2(\Omega)}^2
	\|v^\frac{q}{2}(\cdot,s)\|_{L^\frac{2}{q}(\Omega)}^\frac{4}{q}
	+ \|v^\frac{q}{2}(\cdot,s)\|_{L^\frac{2}{q}(\Omega)}^\frac{2r}{q} \bigg\} ds \\[2mm]
	&\le& c_3 \cdot \Big\{ c_1(T) c_2^2(T) + c_2^r(T) \Big\}
	\qquad \mbox{for all } t\in (0,\htm-\tau),
  \eas
  which implies (\ref{22.1}) with (\ref{22.11}) being valid due to (\ref{22.33}).
\qed
\subsection{Analysis of the functional $\io u^p v^q$ for small positive $p$ and certain $q>0$}
We can now proceed to the following lemma which provides some regularity information that will be fundamental
for our subsequent analysis.
\begin{lem}\label{lem1}
  Let $p\in (0,1)$ be such that $p<\frac{1}{\chi^2}$ and suppose that $q\in (\qm,\qp)$, where
  \be{1.1}
	\qpm:=\frac{1-p}{2} \Big(1\pm\sqrt{1-p\chi^2}\Big).
  \ee
  Then for all $T>0$ there exists $C(T)>0$ such that
  with $\htm:=\min\{T,\tm\}$ and $\tau:=\min\{1,\frac{1}{3}\tm\}$ we have
  \be{1.2}
	\int_t^{t+\tau} \io u^{p-2} v^q u_x^2 \le C(T)
	\qquad \mbox{for all } t\in (0,\htm-\tau),
  \ee
  as well as
  \be{1.3}
	\int_t^{t+\tau} \io u^p v^{q-2} v_x^2 \le C(T)
	\qquad \mbox{for all } t\in (0,\htm-\tau),
  \ee
  and 
  \be{1.33}
	\sup_{T>0} C(T) < \infty
	\qquad \mbox{if either \eqref{H1} or \eqref{H2} hold.}
  \ee
\end{lem}
\proof
  Using that $u$ and $v$ are both positive in $\bom\times (0,\tm)$, on the basis of (\ref{0}) and several integrations
  by parts we compute
  \bea{1.4}
	\frac{d}{dt} \io u^p v^q
	&=& p \io u^{p-1} v^q \cdot \Big\{ u_{xx} - \chi \Big(\frac{u}{v} v_x\Big)_x - uv + B_1 \Big\} 
	+ q\io u^p v^{q-1} \cdot \Big\{ v_{xx} + uv - v + B_2\Big\} \nn\\[1mm]
	&=& p(1-p) \io u^{p-2} v^q u_x^2 
	- pq \io u^{p-1} v^{q-1} u_x v_x \nn\\
	& & - p(1-p)\chi \io u^{p-1} v^{q-1} u_x v_x
	+ pq\chi \io u^p v^{q-2} v_x^2 \nn\\
	& & - p\io u^p v^{q+1} + p \io B_1 u^{p-1} v^q \nn\\
	& & - pq\io u^{p-1} v^{q-1} u_x v_x
	+ q(1-q) \io u^p v^{q-2} v_x^2 \nn\\
	& & + q \io u^{p+1} v^q
	- q\io u^p v^q 
	+ q\io B_2 u^p v^{q-1} \nn\\[1mm]
	&=& p(1-p) \io u^{p-2} v^q u_x^2
	+ q(p\chi+1-q) \io u^p v^{q-2} v_x^2 \nn\\
	& & - p(\chi-p\chi+2q) \io u^{p-1} v^{q-1} u_x v_x \nn\\
	& & -p\io u^p v^{q+1} + p\io B_1 u^{p-1} v^q \nn\\
	& & + q\io u^{p+1} v^q - q\io u^p v^q + q\io B_2 u^p v^{q-1}
	\qquad \mbox{for all } t\in (0,\tm).
  \eea
  Here in order to estimate the third summand on the right, we note that our assumption (\ref{1.1}) on $q$ warrants that
  \bas
	4q^2 - 4(1-p)q + p(1-p)^2 \chi^2 < 0 
  \eas
  and hence
  \bas
	\frac{p(\chi-p\chi+2q)^2}{4(1-p)} - q(p\chi+1-q)
	&=& \frac{1}{4(1-p)} \cdot \Bigg\{ 
	\Big\{ p\chi^2 + p^3\chi^2 + 4pq^2 -2p^2 \chi^2 + 4pq\chi - 4p^2 q\chi\Big\} \\
	& & - \Big\{ 4pq\chi-4p^2 q\chi + 4q - 4pq - 4q^2 + 4pq^2 \Big\} \Bigg\} \\
	&=& \frac{1}{4(1-p)} \cdot \Big\{ 4q^2 - 4(1-p)q + p(1-p)^2 \chi^2 \Big\} \\[2mm]
	&<& 0,
  \eas
  so that it is possible to pick $\eta\in (0,1)$ suitably close to $1$ such that still
  \be{1.5}
	\frac{p(\chi-p\chi+2q)^2}{4(1-p)\eta}
	< q(p\chi+1-q).
  \ee
  Therefore, by Young's inequality we can estimate
  \bas
	& & \hspace*{-20mm}
	p(1-p) \io u^{p-2} v^q u_x^2
	+ q(p\chi+1-q) \io u^p v^{q-2} v_x^2
	- p(\chi-p\chi+2q) \io u^{p-1} v^{q-1} u_x v_x \\
	&\ge&
	p(1-p) \io u^{p-2} v^q u_x^2
	+ q(p\chi+1-q) \io u^p v^{q-2} v_x^2 \\
	& & - \eta p(1-p) \io u^{p-2} v^q u_x^2
	- \frac{p(\chi-p\chi+2q)^2}{4(1-p)\eta} \io u^p v^{q-2} v_x^2  \\
	&=& c_1 \io u^{p-2} v^q u_x^2 + c_2 \io u^p v^{q-2} v_x^2
	\qquad \mbox{for all } t\in (0,\tm),
  \eas
  where $c_1:=(1-\eta)p(1-p)$ is positive due to the fact that $\eta<1$, and where
  \bas
	c_2:=q(p\chi+1-q) - \frac{p(\chi-p\chi+2q)^2}{4(1-p)\eta} >0
  \eas
  thanks to (\ref{1.5}).\abs
  By dropping four nonnegative summands, on integrating (\ref{1.4}) we thus infer that
  \bea{1.6}
	c_1 \int_t^{t+\tau} \io u^{p-2} v^q u_x^2
	+ c_2 \int_t^{t+\tau} \io u^p v^{q-2} v_x^2
	&\le& \io u^p(\cdot,t+\tau) v^q(\cdot,t+\tau) \nn\\
	& & + p \int_t^{t+\tau} \io u^p v^{q+1}
	+ q \int_t^{t+\tau} \io u^p v^q 
  \eea
  for all $t\in (0,\htm-\tau)$.
  Since (\ref{1.1}) particularly requires that
  \be{1.66}
	q<1-p,
  \ee
  we may use the H\"older inequality to see that 
  \bas
	\io u^p v^q \le |\Omega|^{1-p-q} \bigg\{ \io u\bigg\}^p \bigg\{ \io v\bigg\}^q
	\qquad \mbox{for all } t\in (0,\tm),
  \eas
  which in view of Lemma \ref{lem01} and Lemma \ref{lem02} implies that there exists $c_3(T)>0$ such that
  \be{1.7}
	\io u^p(\cdot,t+\tau)v^q(\cdot,t+\tau) + q \int_t^{t+\tau} \io u^p v^q
	\le c_3(T)
	\qquad \mbox{for all } t\in (\htm-\tau),
  \ee
  where
  \be{1.8}
	\sup_{T>0} c_3(T)<\infty
	\qquad \mbox{if either (\ref{H1}) or (\ref{H2}) holds.}
  \ee
  To estimate the second to last summand in (\ref{1.6}), we recall that Lemma \ref{lem01} moreover yields $c_4(T)>0$
  satisfying
  \be{1.9}
	\int_t^{t+\tau} \io uv \le c_4(T)
	\qquad \mbox{for all } t\in (0,\htm-\tau),
  \ee
  where 
  \be{1.10}
	\sup_{T>0} c_4(T)<\infty
	\qquad \mbox{if either (\ref{H1}) or (\ref{H2}) is satisfied.}
  \ee
  Therefore, once again by the H\"older inequality,
  \bea{1.11}
	p\int_t^{t+\tau} \io u^p v^{q+1}
	&=& p \int_t^{t+\tau} \io (uv)^p v^{q+1-p} \nn\\
	&\le& p \bigg\{ \int_t^{t+\tau} \io uv \bigg\}^p 
	\bigg\{ \int_t^{t+\tau} \io v^\frac{q+1-p}{1-p}\bigg\}^{1-p} \nn\\
	&\le& pc_4^p(T)
	\bigg\{ \int_t^{t+\tau} \io v^\frac{q+1-p}{1-p}\bigg\}^{1-p}
	\qquad \mbox{for all } t\in (0,\htm-\tau),
  \eea
  and again by (\ref{1.66}) we see that $\frac{q+1-p}{1-p} < 2 < 3$ and thus Lemma \ref{lem22} becomes applicable 
  to yield $c_5(T)>0$ such that
  \be{1.12}
	\int_t^{t+\tau} \io v^\frac{q+1-p}{1-p} \le c_5(T)
	\qquad \mbox{for all } t\in (0,\htm-\tau),
  \ee
  with
  \be{1.13}
	\sup_{T>0} c_5(T)<\infty
	\qquad \mbox{if either (\ref{H1}) or (\ref{H2}) is valid.}
  \ee
  In summary, (\ref{1.6}), (\ref{1.7}), (\ref{1.11}) and (\ref{1.12}) entail that
  \bas
	c_1 \int_t^{t+\tau} \io u^{p-2} v^q u_x^2
	+ c_2 \int_t^{t+\tau} \io u^p v^{q-2} v_x^2
	\le C(T):=c_3(T) + pc_4^p(T) c_5^{1-p}(T)
	\quad \mbox{for all } t\in (0,\htm-\tau),
  \eas
  where $C(T)$ satisfies (\ref{1.33}) due to (\ref{1.8}), (\ref{1.10}) and (\ref{1.13}).
\qed
\mysection{Global existence. $L^\infty$ bounds for $u$ and $v$ when (\ref{H2}) holds}\label{sec:global}
As a first application of Lemma \ref{lem1}, merely relying on the first inequality (\ref{1.2}) therein and
the pointwise positivity properties of $v$ from Lemma \ref{lem001}
we shall derive a bound for the first solution component in some superquadratic space-time Lebesgue norm.
\begin{lem}\label{lem2}
  Let $p\in (0,1)$ be such that $p<\frac{1}{\chi^2}$.
  Then for all $T>0$ there exists $C(T)>0$ such that
  \be{2.1}
	\int_t^{t+\tau} \io u^{p+2} \le C(T)
	\qquad \mbox{for all } t\in (0,\htm-\tau),
  \ee
  where $\htm:=\min\{T,\tm\}$ and $\tau:=\min\{1,\frac{1}{3}\tm\}$.  Moreover,
  \be{2.2}
	\sup_{T>0} C(T)<\infty
	\qquad \mbox{if \eqref{H2} holds.}
  \ee
\end{lem}
\proof
  We fix any $q\in (\qm,\qp)$, with $\qpm$ taken as in (\ref{1.1}), and invoke Lemma \ref{lem1} and Lemma \ref{lem01}
  to obtain $c_1(T)>0$ and $c_2(T)>0$ such that
  \be{2.3}
	\int_t^{t+\tau} \io u^{p-2} v^q u_x^2 \le c_1(T)
	\qquad \mbox{for all } t\in (0,\htm-\tau)
  \ee
  and
  \be{2.4}
	\io u(\cdot,t) \le c_2(T)
	\qquad \mbox{for all } t\in (0,\htm),
  \ee
  with
  \be{2.5}
	\sup_{T>0} \Big(c_1(T)+c_2(T)\Big) <\infty
	\qquad \mbox{if (\ref{H2}) holds.}
  \ee
  To exploit (\ref{2.3}), we moreover invoke Lemma \ref{lem001} to find $c_3(T)>0$ such that
  \be{2.6}
	v(x,t) \ge c_3(T)
	\qquad \mbox{for all $x\in\Omega$ and } t\in (0,\htm)
  \ee
  with
    \be{2.7}
	\inf_{T>0} c_4(T)>0
	\qquad \mbox{if (\ref{H2}) is valid.}
  \ee
  Therefore, namely, (\ref{2.3}) entails that
  \bas
	\int_t^{t+\tau} \io [(u^\frac{p}{2})_x]^2
	\le c_4(T):=\frac{p^2}{4} \cdot \frac{c_1(T)}{c_3^q(T)}
	\qquad \mbox{for all } t\in (0,\htm-\tau),
  \eas
  and since the Gagliardo-Nirenberg inequality says that with some $c_5>0$ we have
  \bas
	\int_t^{t+\tau} \io u^{p+2}
	&=& \int_t^{t+\tau} \|u^\frac{p}{2}(\cdot,s)\|_{L^\frac{2(p+2)}{p}(\Omega)}^\frac{2(p+2)}{p} ds \\
	&\le& c_5 \int_t^{t+\tau} \bigg\{
	\Big\| (u^\frac{p}{2})_x(\cdot,s)\Big\|_{L^2(\Omega)}^2
	\|u^\frac{p}{2}(\cdot,s)\|_{L^\frac{2}{p}(\Omega)}^\frac{4}{p}
	+ \|u^\frac{p}{2}(\cdot,s)\|_{L^\frac{2}{p}(\Omega)}^\frac{2(p+2)}{p} \bigg\} ds
  \eas
  for all $t\in (0,\htm-\tau)$, by using (\ref{2.4}) we infer that
  \bas
	\int_t^{t+\tau} \io u^{p+2}
	&\le& c_5 c_2^2(T) \int_t^{t+\tau} \io (u^\frac{p}{2})_x^2 + c_5 c_2^{p+2}(T) \\
	&\le& c_5 c_2^2(T) c_4(T) + c_5 c_2^{p+2}(T)
	\qquad \mbox{for all } t\in (0,\htm-\tau).
  \eas
  Combined with (\ref{2.5}) and (\ref{2.7}) this establishes (\ref{2.1}) and (\ref{2.2}).
\qed
In the considered spatially one-dimensional case, the latter property turns out to be sufficient for the derivation
of bounds for $v_x$ in $L^r(\Omega)$ for suitably small $r>1$.
\begin{lem}\label{lem3}
  Let $r\in (1,\frac{3}{2})$ be such that $r<1+\frac{1}{2\chi^2}$.
  Then for all $T>0$ there exists $C(T)>0$ such that with
  $\htm:=\min\{T,\tm\}$ we have
  \be{3.1}
	\|v_x(\cdot,t)\|_{L^r(\Omega)} \le C(T)
	\qquad \mbox{for all } t\in (0,\htm),
  \ee
  where
  \be{3.2}
	\sup_{T>0} C(T)<\infty
	\qquad \mbox{if \eqref{H2} holds.}
  \ee
\end{lem}
\proof
  Once more writing $\tau:=\min\{1,\frac{1}{3}\tm\}$, from Lemma \ref{lem_loc} we know that
  \bas
	c_1:=\sup_{t\in (0,\tau]} \|v_x(\cdot,t)\|_{L^r(\Omega)}
  \eas
  is finite, whence for estimating
  \bas
	M(T'):=\sup_{t\in (0,T')} \|v_x(\cdot,t)\|_{L^r(\Omega)}
	\qquad \mbox{for } T'\in (\tau,\htm)
  \eas
  it will be sufficient to derive appropriate bounds of $v_x(\cdot,t)$ in $L^r(\Omega)$ for $t\in (\tau,T')$ only.
  To this end, given any such $t$ we represent $v_x(\cdot,t)$ according to
  \bea{3.3}
	v_x(\cdot,t)
	= \partial_x e^{\tau(\Delta-1)} v(\cdot,t-\tau)
	+ \int_{t-\tau}^t \partial_x e^{(t-s)(\Delta-1)} u(\cdot,s)v(\cdot,s) ds
	+ \int_{t-\tau}^t \partial_x e^{(t-s)(\Delta-1)} B_2(\cdot,s) ds
  \eea
  and recall that due to known smoothing properties of the Neumann heat semigroup (\cite{win_JDE}) we can find
  $c_2>0$ such that for all $\varphi\in C^0(\bom)$,
  \be{3.4}
	\|\partial_x e^{\sigma\Delta}\varphi\|_{L^r(\Omega)} 
	\le c_2 \sigma^{-\frac{1}{2}-\frac{1}{2}(1-\frac{1}{r})} \|\varphi\|_{L^1(\Omega)}
	\qquad \mbox{for all } \sigma\in (0,1).
  \ee
  Therefore,
  \bea{3.5}
	\Big\|\partial_x e^{\tau(\Delta-1)} v(\cdot,t-\tau)\Big\|_{L^r(\Omega)}
	&\le& c_2 e^{-\tau} \cdot \tau^{-\frac{1}{2}-\frac{1}{2}(1-\frac{1}{r})} \|v(\cdot,t-\tau)\|_{L^1(\Omega)} \nn\\
	&\le& c_2 c_3(T) \tau^{-\frac{1}{2}-\frac{1}{2}(1-\frac{1}{r})},
  \eea
  where $c_3(T)>0$ has been chosen in such a way that in accordance with Lemma \ref{lem02} we have
  \be{3.6}
	\|v(\cdot,s)\|_{L^1(\Omega)} \le c_3(T)
	\qquad \mbox{for all } s\in (0,\htm),
  \ee
  and such that
  \be{3.7}
	\sup_{T>0} c_3(T)<\infty
	\qquad \mbox{if (\ref{H2}) holds.}
  \ee
  Next, again by (\ref{3.4}),
  \bea{3.8}
	\bigg\| \int_{t-\tau}^t \partial_x e^{(t-s)(\Delta-1)} B_2(\cdot,s) ds \bigg\|_{L^r(\Omega)} ds
	&\le& c_2 \int_{t-\tau}^t (t-s)^{-\frac{1}{2}-\frac{1}{2}(1-\frac{1}{r})} \|B_2(\cdot,s)\|_{L^1(\Omega)} ds \nn\\
	&\le& c_2|\Omega| \|B_2\|_{L^\infty(\Omega\times (0,\infty))} 
	\int_0^\tau \sigma^{-\frac{1}{2}-\frac{1}{2}(1-\frac{1}{r})} d\sigma \nn\\[2mm]
	&=& c_2|\Omega| \|B_2\|_{L^\infty(\Omega\times (0,\infty))} 
	\cdot 2r \tau^\frac{1}{2r}
  \eea
  as well as
  \be{3.9}
	\bigg\| \int_{t-\tau}^t \partial_x e^{(t-s)(\Delta-1)} u(\cdot,s) v(\cdot,s) ds \bigg\|_{L^r(\Omega)} ds
	\le c_2 \int_{t-\tau}^t (t-s)^{-\frac{1}{2}-\frac{1}{2}(1-\frac{1}{r})} 
	\|u(\cdot,s)v(\cdot,s)\|_{L^1(\Omega)} ds.
  \ee
  In order to further estimate the latter integral, we make use of our restrictions $r<\frac{3}{2}$ and
  $r<1+\frac{1}{2\chi^2}$ which enable us to pick some $p\in (0,1)$ satisfying $p<\frac{1}{\chi^2}$ and $p>2(r-1)$.
  Then by means of the H\"older inequality we see that
  \be{3.10}
	\|u(\cdot,s)v(\cdot,s)\|_{L^1(\Omega)}
	\le \|u(\cdot,s)\|_{L^{p+2}(\Omega)} \|v(\cdot,s)\|_{L^\frac{p+2}{p+1}(\Omega)}
	\qquad \mbox{for all } s\in (0,\tm),
  \ee
  where the Gagliardo-Nirenberg inequality provides $c_4>0$ and $a\in (0,1)$ fulfilling
  \bas
	\|v(\cdot,s)\|_{L^\frac{p+2}{p+1}(\Omega)}
	\le c_4\|v_x(\cdot,s)\|_{L^r(\Omega)}^a \|v(\cdot,s)\|_{L^1(\Omega)}^{1-a}
	+ c_4\|v(\cdot,s)\|_{L^1(\Omega)}
	\qquad \mbox{for all } s\in (0,\tm).
  \eas
  In light of (\ref{3.6}) and the definition of $M(T')$, from (\ref{3.10}) and (\ref{3.9}) we thus obtain that
  \bas
	\|u(\cdot,s)v(\cdot,s)\|_{L^1(\Omega)}
	\le \|u(\cdot,s)\|_{L^{p+2}(\Omega)} \cdot 
	\Big\{ c_4 c_3^{1-a}(T) M^a(T') + c_4 c_3(T)\Big\},
  \eas
  so that once again invoking the H\"older inequality we infer that
  \bea{3.11}
	& & \hspace*{-20mm}
	\bigg\| \int_{t-\tau}^t \partial_x e^{(t-s)(\Delta-1)} u(\cdot,s) v(\cdot,s) ds \bigg\|_{L^r(\Omega)} ds \nn\\
	&\le& c_2 \cdot \Big\{ c_4 c_3^{1-a}(T) M^a(T') + c_4 c_3(T)\Big\}
	\cdot \int_{t-\tau}^t (t-s)^{-\frac{1}{2}-\frac{1}{2}(1-\frac{1}{r})} \|u(\cdot,s)\|_{L^{p+2}(\Omega)} ds \nn\\
	&\le& c_2 \cdot \Big\{ c_4 c_3^{1-a}(T) M^a(T') + c_4 c_3(T)\Big\}
	\cdot \bigg\{ 
	\int_{t-\tau}^t (t-s)^{-[\frac{1}{2}+\frac{1}{2}(1-\frac{1}{r})] \cdot \frac{p+2}{p+1}} ds
	\bigg\}^\frac{p+1}{p+2} \times \nn\\
	& & \hspace*{58mm}
	\times \bigg\{ \int_{t-\tau}^t \|u(\cdot,s)\|_{L^{p+2}(\Omega)}^{p+2} ds \bigg\}^\frac{1}{p+2}.
  \eea
  Since herein our assumption $p>2(r-1)$ warrants that
  \bas
	\Big[\frac{1}{2}+\frac{1}{2}\Big(1-\frac{1}{r}\Big)\Big] \cdot \frac{p+2}{p+1}
	= \frac{2r-1}{2r} \cdot \Big(1+\frac{1}{p+1}\Big)
	< \frac{2r-1}{2r} \cdot \Big(1+\frac{1}{2r-1}\Big)
	= 1,
  \eas
  and since from Lemma \ref{lem2} we know that
  \bas
	\int_{t-\tau}^t \|u(\cdot,s)\|_{L^{p+2}(\Omega)}^{p+2} ds
	\le c_5(T)
  \eas
  with some $c_5(T)>0$ satisfying
  \bas
	\sup_{T>0} c_5(T)<\infty
	\qquad \mbox{if (\ref{H2}) holds,}
  \eas
  it follows from (\ref{3.11}) that with a certain $c_6(T)>0$ we have
  \bas
	\bigg\| \int_{t-\tau}^t \partial_x e^{(t-s)(\Delta-1)} u(\cdot,s) v(\cdot,s) ds \bigg\|_{L^r(\Omega)} ds
	\le c_6(T) \cdot \Big\{ M^a(T')+1\Big\},
  \eas
  where
  \be{3.13}
	\sup_{T>0} c_6(T)<\infty
	\qquad \mbox{if (\ref{H2}) is valid.}
  \ee
  Together with (\ref{3.5}), (\ref{3.8}) and (\ref{3.3}), this shows that
  \be{3.14}
	\|v_x(\cdot,t)\|_{L^r(\Omega)} \le c_7(T) \cdot \Big\{ M^a(T')+1 \Big\}
	\qquad \mbox{for all } t\in (\tau,T')
  \ee
  with some $c_7(T)>0$ which due to (\ref{3.7}) and (\ref{3.13}) is such that
  \be{3.15}
	\sup_{T>0} c_7(T)<\infty
	\qquad \mbox{if (\ref{H2}) holds.}
  \ee
  In view of our definition of $c_1$, (\ref{3.14}) entails that if we let
  $c_8(T):=\max\{c_7(T),1\}$ then
  \bas
	M(T') \le c_8(T) \cdot \Big\{ M^a(T')+1\Big\}
	\qquad \mbox{for all } T'\in (\tau,\htm)
  \eas
  and thus, since $a<1$,
  \bas
	M(T') \le \max \Big\{ 1 \, , \, (2c_8(T))^\frac{1}{1-a}\Big\}
	\qquad \mbox{for all } T'\in (\tau,\htm).
  \eas
  Combined with (\ref{3.15}), this establishes (\ref{3.1}) and (\ref{3.2}).
\qed
Now, the latter provides sufficient regularity of the inhomogeneity $h$ appearing in the identity $u_t=u_{xx}+h$
in (\ref{0}), that is, of $h:=-\chi(\frac{u}{v}v_x)_x -uv+B_1$, and especially in the crucial cross-diffusive 
first summand therein. This is obtained by the following statement which beyond boundedness of $u$, as required
for extending the solution via Lemma \ref{lem_loc}, moreover asserts a favorable equicontinuity feature
of $u$ that will be useful in verifying the uniform decay property claimed in Theorem \ref{theo51}.
\begin{lem}\label{lem4}
  Let $\gamma\in (0,\frac{1}{3})$ be such that $\gamma<\frac{1}{1+2\chi^2}$.
  Then for all $T>0$ there exists $C(T)>0$ with the properties that with 
  $\htm:=\min\{T,\tm\}$ and $\tau:=\min\{1,\frac{1}{3}\tm\}$ we have
  \be{4.1}
	\|u(\cdot,t)\|_{C^\gamma(\bom)} \le C(T),
	\qquad \mbox{for all } t\in (\tau,\htm)
  \ee
  and
  \be{4.2}
	\sup_{T>0} C(T)<\infty
	\qquad \mbox{if \eqref{H2} holds.}
  \ee
\end{lem}
\proof
  Since $\gamma<\frac{1}{3}$ and $\gamma<\frac{1}{1+2\chi^2}$, it is possible to fix $r>1$ such that $r<\frac{3}{2}$
  and $r<1+\frac{1}{2\chi^2}$, and such that $1-\frac{1}{r}>\gamma$.
  This enables us to choose some $\alpha\in (0,\frac{1}{2})$ sufficiently close to $\frac{1}{2}$ such that still
  $2\alpha-\frac{1}{r}>\gamma$, which in turn ensures that the sectorial realization of $A:=-(\cdot)_{xx}+1$
  under homogeneous Neumann boundary conditions in $L^r(\Omega)$ has the domain of its fractional power $A^\alpha$
  satisfy $D(A^\alpha) \hra C^\gamma(\bom)$ (\cite{henry}), meaning that
  \be{4.3}
	\|\varphi\|_{C^\gamma(\bom)} \le c_1 \|A^\alpha \varphi\|_{L^r(\Omega)}
	\qquad \mbox{for all } \varphi \in C^1(\bom)
  \ee
  with some $c_1>0$.
  Moreover, combining known regularization estimates for the associated semigroup
  $(e^{-tA})_{t\ge 0}\equiv (e^{-t} e^{t\Delta})_{t\ge 0}$ (\cite{friedman}, \cite{win_JDE}) we can find positive 
  constants $c_2$ and $c_3$ such that for all $t\in (0,1)$ we have
  \be{4.4}
	\|A^\alpha e^{-tA} \varphi\|_{L^r(\Omega)}
	\le c_2 t^{-\alpha-\frac{1}{2}(1-\frac{1}{r})} \|\varphi\|_{L^1(\Omega)}
	\qquad \mbox{for all } \varphi\in C^0(\bom)
  \ee
  and
  \be{4.5}
	\|A^\alpha e^{-tA}\varphi_x\|_{L^r(\Omega)} 
	\le c_3 t^{-\alpha-\frac{1}{2}} \|\varphi\|_{L^r(\Omega)}
	\qquad \mbox{for all } \varphi\in C^1(\bom)
	\mbox{ such that $\varphi_x=0$ on } \pO.
  \ee
  Now to estimate
  \bas
	M(T'):=\sup_{t\in (\tau,T')} \|u(\cdot,t)\|_{C^\gamma(\bom)}
	\qquad \mbox{for } T' \in (\tau,\htm),
  \eas
  we use a variation-of-constants representation associated with the identity
  \bas
	u_t=-Au - \chi \nabla\cdot \Big(\frac{u}{v} v_x\Big)_x - uv + b_1(x,t) + u,
	\qquad x\in\Omega, \ t\in (0,\tm),
  \eas
  to see that thanks to (\ref{4.3}),
  \bea{4.6}
	\frac{1}{c_1} \|u(\cdot,t)\|_{C^\gamma(\bom)}
	&\le& \|A^\alpha u(\cdot,t)\|_{L^r(\Omega)} \nn\\
	&\le& \Big\| A^\alpha e^{-\tau A} u(\cdot,t-\tau)\Big\|_{L^r(\Omega)}
	+ \chi \int_{t-\tau}^t \Big\| A^\alpha e^{-(t-s)A} \Big(\frac{u(\cdot,s)}{v(\cdot,s)} v_x(\cdot,s)\Big)_x
		\Big\|_{L^r(\Omega)} ds \nn\\
	& & + \int_{t-\tau}^t \Big\| A^\alpha e^{-(t-s)A} u(\cdot,s) v(\cdot,s) \Big\|_{L^r(\Omega)} ds 
	+ \int_{t-\tau}^t \Big\| A^\alpha e^{-(t-s)A} B_1(\cdot,s)\Big\|_{L^r(\Omega)} ds \nn\\
	& & + \int_{t-\tau}^t \Big\| A^\alpha e^{-(t-s)A} u(\cdot,s)\Big\|_{L^r(\Omega)} ds
	\qquad \mbox{for all } t\in (2\tau,\tm).
  \eea
  Here by (\ref{4.4}) we see that
  \bea{4.7}
	\|A^\alpha e^{-\tau A} u(\cdot,t-\tau)\|_{L^r(\Omega)}
	&\le& c_2 \tau^{-\alpha-\frac{1}{2}(1-\frac{1}{r})} \|u(\cdot,t-\tau)\|_{L^1(\Omega)} \nn\\
	&\le& c_2 c_4(T) \tau^{-\alpha-\frac{1}{2}(1-\frac{1}{r})}
	\qquad \mbox{for all } t\in (2\tau,\htm),
  \eea
  where according to Lemma \ref{lem01} we have taken $c_4(T)>0$ such that
  \be{4.8}
	\|u(\cdot,t)\|_{L^1(\Omega)} \le c_4(T)
	\qquad \mbox{for all } t\in (0,\htm)
  \ee
  and that
  \be{4.9}
	\sup_{T>0} c_4(T)<\infty
	\qquad \mbox{if (\ref{H2}) holds.}
  \ee
  Moreover, in view of our restrictions on $r$ we see that Lemma \ref{lem3} applies so as to yield $c_5(T)>0$ satisfying
  \be{4.10}
	\|v_x(\cdot,t)\|_{L^r(\Omega)} \le c_5(T)
	\qquad \mbox{for all } t\in (0,\htm)
  \ee
  and
  \be{4.11}
	\sup_{T>0} c_5(T)<\infty
	\qquad \mbox{if (\ref{H2}) is valid,}
  \ee
  which combined with the outcome of Lemma \ref{lem02} and the continuity of the embedding $W^{1,r}(\Omega) \hra
  L^\infty(\Omega)$ shows that there exists $c_6(T)>0$ such that
  \be{4.12}
	\|v(\cdot,t)\|_{L^\infty(\Omega)} \le c_6(T)
	\qquad \mbox{for all } t\in (0,\htm)
  \ee
  with
  \be{4.13}
	\sup_{T>0} c_6(T)<\infty
	\qquad \mbox{if (\ref{H2}) holds.}
  \ee
  Therefore, in the third integral on the right of (\ref{4.6}) we may use (\ref{4.4}) and again (\ref{4.8}) to estimate
  \bea{4.14}
	\int_{t-\tau}^t \Big\|A^\alpha e^{-(t-s)A} u(\cdot,s)v(\cdot,s)\Big\|_{L^r(\Omega)} ds
	&\le& c_2 \int_{t-\tau}^t (t-s)^{-\alpha-\frac{1}{2}(1-\frac{1}{r})}
	\|u(\cdot,s)v(\cdot,s)\|_{L^1(\Omega)} ds \nn\\
	&\le& c_2 \int_{t-\tau}^t (t-s)^{-\alpha-\frac{1}{2}(1-\frac{1}{r})}
	\|u(\cdot,s)\|_{L^1(\Omega)} \|v(\cdot,s)\|_{L^\infty(\Omega)} ds \nn\\
	&\le& c_2 c_4(T) c_6(T) \int_{t-\tau}^t (t-s)^{-\alpha-\frac{1}{2}(1-\frac{1}{r})} ds \nn\\
	&=& c_2 c_4(T) c_6(T) c_7
	\qquad \mbox{for all } t\in (2\tau,\htm),
  \eea
  with $c_7:=\int_0^\tau \sigma^{-\alpha-\frac{1}{2}(1-\frac{1}{r})}d\sigma$ being finite since clearly
  $\alpha+\frac{1}{2}(1-\frac{1}{r})<\alpha+\frac{1}{2}<1$.\\
  Likewise, upon two further applications of (\ref{4.4}) we obtain from the boundedness of $B_1$ and (\ref{4.8}) that
  \bea{4.15}
	\int_{t-\tau}^t \Big\| A^\alpha e^{-(t-s)A} B_1(\cdot,s)\Big\|_{L^r(\Omega)} ds
	&\le& c_2 \int_{t-\tau}^t (t-s)^{-\alpha-\frac{1}{2}(1-\frac{1}{r})} \|B_1(\cdot,s)\|_{L^1(\Omega)} ds \nn\\
	&\le& c_2 |\Omega| \|B_1\|_{L^\infty(\Omega\times (0,\infty))} 
	\int_{t-\tau}^t (t-s)^{-\alpha-\frac{1}{2}(1-\frac{1}{r})} ds \nn\\
	&=& c_2 |\Omega| \|B_1\|_{L^\infty(\Omega\times (0,\infty))} \cdot c_7
	\qquad \mbox{for all } t\in (2\tau,\htm)
  \eea
  and that
  \bea{4.16}
	\int_{t-\tau}^t \Big\| A^\alpha e^{-(t-s)A} u(\cdot,s)\Big\|_{L^r(\Omega)} ds
	&\le& c_2 \int_{t-\tau}^t (t-s)^{-\alpha-\frac{1}{2}(1-\frac{1}{r})} \|u(\cdot,s)\|_{L^1(\Omega)} ds \nn\\
	&\le& c_2 c_4(T) c_7
	\qquad \mbox{for all } t\in (2\tau,\htm).
  \eea
  Finally, in the second summand on the right-hand side in (\ref{4.6}) we use that due to Lemma \ref{lem001},
  \bas
	v(x,t) \ge c_8(T) 
	\qquad \mbox{for all $x\in\Omega$ and } t\in (0,\htm)
  \eas
  with some $c_8(T)>0$ fulfilling
  \be{4.17}
	\inf_{T>0} c_8(T)>0
	\qquad \mbox{if (\ref{H2}) holds.}
  \ee
  From (\ref{4.5}) and (\ref{4.10}) we therefore obtain that
  \bea{4.18}
	& & \hspace*{-20mm}
	\chi \int_{t-\tau}^t \Big\| A^\alpha e^{-(t-s)A} \Big(\frac{u(\cdot,s)}{v(\cdot,s)} v_x(\cdot,s)\Big)_x
		\Big\|_{L^r(\Omega)} ds \nn\\
	&\le& \chi c_3 \int_{t-\tau}^t (t-s)^{-\alpha-\frac{1}{2}} \Big\|\frac{u(\cdot,s)}{v(\cdot,s)} v_x(\cdot,s)
		\Big\|_{L^r(\Omega)} ds \nn\\
	&\le& \chi c_3 \int_{t-\tau}^t (t-s)^{-\alpha-\frac{1}{2}} 
	\|u(\cdot,s)\|_{L^\infty(\Omega)} \Big\|\frac{1}{v(\cdot,s)}\Big\|_{L^\infty(\Omega)}
	\|v_x(\cdot,s)\|_{L^r(\Omega)} ds \nn\\
	&\le& \frac{\chi c_3 c_5(T)}{c_8(T)} \|u\|_{L^\infty(\Omega\times (\tau,T'))}
	\int_{t-\tau}^t (t-s)^{-\alpha-\frac{1}{2}} ds \nn\\
	&=& \frac{\chi c_3 c_5(T)}{c_8(T)} \|u\|_{L^\infty(\Omega\times (\tau,T'))}
	\cdot \frac{\tau^{\frac{1}{2}-\alpha}}{\frac{1}{2}-\alpha}
	\qquad \mbox{for all } t\in (2\tau,T').
  \eea
  In conclusion, (\ref{4.7}), (\ref{4.14}), (\ref{4.15}), (\ref{4.16}) and (\ref{4.18}) show that (\ref{4.6})
  leads to the inequality
  \be{4.188}
	\|u(\cdot,t)\|_{C^\gamma(\bom)}
	\le c_9(T) \|u\|_{L^\infty(\Omega\times (\tau,T'))} + c_9(T)
	\qquad \mbox{for all } t\in (2\tau,T')
  \ee 
  with some $c_9(T)>0$ about which due to (\ref{4.9}), (\ref{4.11}), (ref{4.13}) and (\ref{4.17}) we know that
  \be{4.19}
	\sup_{T>0} c_9(T)<\infty
	\qquad \mbox{if (\ref{H2}) holds.}
  \ee
  Now, by compactness of the first in the embeddings $C^\gamma(\bom)\hra L^\infty(\Omega) \hra L^1(\Omega)$,
  according to an associated Ehrling lemma it is possible to pick $c_{10}(T)>0$ such that
  \bas
	\|\varphi\|_{L^\infty(\Omega)}
	\le \frac{1}{2c_9(T)} \|\varphi\|_{C^\gamma(\bom)}
	+ c_{10}(T) \|\varphi\|_{L^1(\Omega)}
	\qquad \mbox{for all } \varphi\in C^\gamma(\bom),
  \eas
  where thanks to (\ref{4.19}) it can clearly be achieved that
  \be{4.20}
	\sup_{T>0} c_{10}(T)<\infty,
	\qquad \mbox{provided that (\ref{H2}) holds.}
  \ee
  Therefore, (\ref{4.188}) together with (\ref{4.8}) implies that
  \bas
	\|u(\cdot,t)\|_{C^\gamma(\bom)}
	&\le& \frac{1}{2} \sup_{s\in (\tau,T')} \|u(\cdot,s)\|_{C^\gamma(\bom)}
	+ c_{10}(T) c_4(T) + c_9(T) \\
	&\le& \frac{1}{2} M(T')
	+ c_{10}(T) c_4(T) + c_9(T) 
	\qquad \mbox{for all } t\in (2\tau,T')
  \eas
  and that hence with $c_{11}:=\sup_{t\in (\tau,2\tau]} \|u(\cdot,t)\|_{C^\gamma(\bom)}$ we have
  \bas
	M(T')
	&\le& c_{11} + \sup_{t\in (2\tau,T')} \|u(\cdot,t)\|_{C^\gamma(\bom)} \\
	&\le& c_{11} + \frac{1}{2} M(T') + c_{10}(T) c_4(T) + c_9(T).
  \eas
  Thus, 
  \bas
	M(T') \le 2\cdot \Big( c_{11} + c_{10}(T) c_4(T) + c_9(T)\Big)
	\qquad \mbox{for all } T'\in (\tau,\htm),
  \eas
  which on letting $T'\nearrow \htm$ yields (\ref{4.1}) with some
  $C(T)>0$ satisfying (\ref{4.2}) because of (\ref{4.20}), (\ref{4.8}) and (\ref{4.19}).
\qed
\subsection{Proof of Theorem \ref{theo5}}
By collecting the above positivity and regularity information, we immediately obtain global extensibility
of our local-in-time solution:\abs
\proofc of Theorem \ref{theo5}. \quad
  Combining Lemma \ref{lem4} with Lemma \ref{lem001} and Lemma \ref{lem3} shows that in (\ref{ext_crit}), the second
  alternative cannot occur, so that actually $\tm=\infty$ and hence all statements result from Lemma \ref{lem_loc}.
\qed
\subsection{Proof of Theorem \ref{theo51}}
In view of the above statements on independence of all essential estimates from $T$ when (\ref{H2}) holds,
for the verification of the qualitative properties in Theorem \ref{theo51} only one further ingredient
is needed which can be obtained by a refined variant of an argument from Lemma \ref{lem01}.
\begin{lem}\label{lem33}
  If \eqref{H2} and \eqref{H1'} are satisfied, then
  \be{33.1}
	\io u(\cdot,t) \to 0
	\qquad \mbox{as } t\to\infty.
  \ee
\end{lem}
\proof
  Since Lemma \ref{lem001} provides $c_1>0$ such that $v\ge c_1$ in $\Omega\times (0,\infty)$, once more integrating the
  first equation in (\ref{0}) we obtain that
  \bas
	\frac{d}{dt} \io u = - \io uv + \io B_1
	\le -c_1 \io u + \io B_1
	\qquad \mbox{for all } t>0.
  \eas
  In view of the hypothesis (\ref{H1'}), the claim therefore results by an application of Lemma \ref{lem_espejo_win}.
\qed
We can thereby prove our main result on large time behavior in (\ref{0}), (\ref{0b}), (\ref{0i}) in presence of the
hypothesis (\ref{H2}).\abs
\proofc of Theorem \ref{theo51}. \quad
   Assuming that (\ref{H2}) be valid, from Lemma \ref{lem4} we obtain $\gamma>0$ and $c_1>0$ such that
  \be{51.5}
	\|u(\cdot,t)\|_{C^\gamma(\bom)} \le c_1
	\qquad \mbox{for all } t>1.
  \ee
  This immediately implies (\ref{51.1}), whereas
  the inequalities in (\ref{51.2}) result from Lemma \ref{lem001}, Lemma \ref{lem02} and Lemma \ref{lem3},
  again because $W^{1,r}(\Omega) \hra L^\infty(\Omega)$ for arbitrary $r>1$.
  Finally, as the Arzel\`a-Ascoli theorem says that (\ref{51.5}) implies precompactness
  of $(u(\cdot,t))_{t>1}$ in $L^\infty(\Omega)$, the outcome of Lemma \ref{lem33}, asserting that (\ref{H1'}) entails
  decay of $u(\cdot,t)$ in $L^1(\Omega)$ as $t\to\infty$, actually means that we must even have
  $u(\cdot,t)\to 0$ in $L^\infty(\Omega)$ as $t\to\infty$ in this case.
\qed
\mysection{Bounds for $v$ under the assumption (\ref{H1}). Proof of Theorem \ref{theo52}} \label{sec:longtime}
In order to prove Theorem \ref{theo52}, we evidently may no longer rely on any global positivity property of
$v$, which in view of the singular taxis term in (\ref{0}) apparently reduces our information on regularity
of $u$ to a substantial extent. 
Our approach will therefore alternatively focus on the derivation of further bounds for $v$ by merely using the
second equation in (\ref{0}) together with the class of fundamental estimates from Lemma \ref{lem1}, 
taking essential advantage from the freedom to choose the parameters $p$ and $q$ there within a suitably large range.\abs
Our argument will at its core be quite simple in that it is built on a straightforward $L^r$ testing procedure
(see Lemma \ref{lem13}); however, for adequately estimating the crucial integrals $\io uv^r$ appearing therein 
we will create an iterative setup which allows the eventual choosing of an arbitrarily large $r$ whenever
$\chi$ satisfies the smallness condition from Theorem \ref{theo52}.\abs
Let us first reformulate the outcome of Lemma \ref{lem1} in a version convenient for our purpose.
\begin{lem}\label{lem10}
  Assume that \eqref{H1} holds, and let $p\in (0,1)$ and $q>0$ be such that $p<\frac{1}{\chi^2}$
  and $q\in (\qm,\qp)$ with $\qpm$ as given by \eqref{1.1}.
  Then there exists $C>0$ such that
  \be{10.1}
	\int_t^{t+1} \io \left[\Big(u^\frac{p}{2} v^\frac{q}{2}\Big)_x\right]^2 \le C
	\qquad \mbox{for all } t>0.
  \ee
\end{lem}
\proof
  Since
  \bas
	\left[\Big(u^\frac{p}{2} v^\frac{q}{2}\Big)_x\right]^2
	&=& \Big(\frac{p}{2} u^\frac{p-2}{2} v^\frac{q}{2} u_x + \frac{q}{2} u^\frac{p}{2} v^\frac{q-2}{2} v_x\Big)^2 \\
	&\le& \frac{p^2}{2} u^{p-2} v^q u_x^2 + \frac{q^2}{2} u^p v^{q-2} v_x^2
	\qquad \mbox{in } \Omega\times (0,\infty),
  \eas
  by Young's inequality, this is an immediate consequence of Lemma \ref{lem1}.
\qed
A zero-order estimate for the coupled quantities appearing in the preceding lemma can be achieved
by combining Lemma \ref{lem01} with a supposedly known bound for $v$ in $L^{\rs}(\Omega)$ 
in a straightforward manner.
\begin{lem}\label{lem11}
  Assume that (\ref{H1}) holds, and let $\rs\ge 1, p>0$ and $q>0$. Then there exists $C>0$ with the property that if
  with some $K>0$ we have
  \be{11.1}
	\|v(\cdot,t)\|_{L^{\rs}(\Omega)} \le K
	\qquad \mbox{for all } t>0,
  \ee
  then
  \be{11.2}
	\Big\|u^\frac{p}{2}(\cdot,t)v^\frac{q}{2}(\cdot,t)\Big\|_{L^\frac{2\rs}{p\rs+q}(\Omega)} 
	\le CK^\frac{q}{2}
	\qquad \mbox{for all } t>0.
  \ee
\end{lem}
\proof
  According to the hypothesis (\ref{H1}), from Lemma \ref{lem01} we know that
  \bas
	\|u(\cdot,t)\|_{L^1(\Omega)} \le c_1
	\qquad \mbox{for all } t>0
  \eas
  with some $c_1>0$. By the H\"older inequality we therefore obtain that
  \bas
	\big\|u^\frac{p}{2}v^\frac{q}{2}\big\|_{L^\frac{2\rs}{p\rs+q}(\Omega)}
	=  \bigg\{ \io u^\frac{p\rs}{p\rs+q} v^\frac{q\rs}{p\rs+q} \bigg\}^\frac{p\rs+q}{2\rs}
	\le \bigg\{ \io u\bigg\}^\frac{p}{2} \cdot \bigg\{ \io v^{\rs} \bigg\}^\frac{q}{2} 
	\le c_1^\frac{p}{2} K^\frac{q}{2}
	\qquad \mbox{for all } t>0
  \eas
  due to (\ref{11.1}).
\qed
We can thereby achieve the following estimate for the crucial term $\io uv^r$ appearing in Lemma \ref{lem13}
below, for certain $r$ depending on the invested integrability parameter $\rs$.
\begin{lem}\label{lem12}
  Assume \eqref{H1} and suppose that there exists $\rs\ge 1$ such that
  \be{12.01}
	\sup_{t>0} \|v(\cdot,t)\|_{L^{\rs}(\Omega)} < \infty.
  \ee
  Moreover, let $p\in (0,1)$ be such that $p<\frac{1}{\chi^2}$, and with $\qpm$ as given by (\ref{1.1}), let
  $q\in (\qm,\qp)$ satisfy
  \be{12.1}
	q \le \frac{p(p+1)}{1-p} \cdot \rs.
  \ee
  Then there exists $C>0$ such that
  \be{12.2}
	\int_t^{t+1} \io uv^\frac{q}{p} \le C
	\qquad \mbox{for all } t>0.
  \ee
\end{lem}
\proof
  From Lemma \ref{lem11} we know that due to (\ref{12.01}) we can pick $c_1>0$ such that
  \be{12.3}
	\Big\|u^\frac{p}{2}(\cdot,t) v^\frac{q}{2}(\cdot,t)\Big\|_{L^\frac{2\rs}{p\rs+q}(\Omega)} \le c_1
	\qquad \mbox{for all } t>0,
  \ee
  and since (\ref{12.1}) warrants that
  \bas
	\frac{2q}{p[(p+1)\rs+q]}
	= \frac{2}{p\cdot \Big[\frac{(p+1)\rs}{q}+1\Big]}
	\le \frac{2}{p\cdot\Big[\frac{1-p}{p}+1\Big]}
	=2,
  \eas
  we may combine the outcome of Lemma \ref{lem10} with Young's inequality to obtain $c_2>0$ fulfilling
  \be{12.4}
	\int_t^{t+1} \Big\| \Big( u^\frac{p}{2}(\cdot,s) v^\frac{q}{2}(\cdot,s)\Big)_x
		\Big\|_{L^2(\Omega)}^\frac{2q}{p[(p+1)\rs+q]} ds
	\le c_2 
	\qquad \mbox{for all } t>0.
  \ee
  As the Gagliardo-Nirenberg inequality provides $c_3>0$ such that
  \bas
	\|\varphi\|_{L^\frac{2}{p}(\Omega)}^\frac{2}{p}
	\le c_3\|\varphi_x\|_{L^2(\Omega)}^\frac{2q}{p[(p+1)\rs+q]}
	\|\varphi\|_{L^\frac{2\rs}{p\rs+q}(\Omega)}^\frac{2(p+1)\rs}{p[(p+1)\rs+q]}
	+ c_3 \|\varphi\|_{L^\frac{2\rs}{p\rs+q}(\Omega)}^\frac{2}{p}
	\qquad \mbox{for all } \varphi\in W^{1,2}(\Omega),
  \eas
  combining (\ref{12.4}) with (\ref{12.3}) we thus infer that
  \bas
	\int_t^{t+1} \io uv^\frac{q}{p}
	&=& \int_t^{t+1} \Big\|u^\frac{p}{2}(\cdot,s) v^\frac{q}{2}(\cdot,s)\Big\|_{L^\frac{2}{p}(\Omega)}^\frac{2}{p} ds
		\\
	&\le& c_3 \int_t^{t+1} \Big\| \Big(u^\frac{p}{2}(\cdot,s)v^\frac{q}{2}(\cdot,s)\Big)_x
		\Big\|_{L^2(\Omega)}^\frac{2q}{p[(p+1)\rs+q]}
	\Big\| u^\frac{p}{2}(\cdot,s)v^\frac{q}{2}(\cdot,s)
		\Big\|_{L^\frac{2\rs}{p\rs+q}(\Omega)}^\frac{2(p+1)\rs}{p[(p+1)\rs+q]} ds \\
	& & + c_3 \int_t^{t+1} \Big\|u^\frac{p}{2}(\cdot,s)v^\frac{q}{2}(\cdot,s)
		\Big\|_{L^\frac{2\rs}{p\rs+q}(\Omega)}^\frac{2}{p} ds \\[2mm]
	&\le& c_3 \cdot c_2 c_1^\frac{2(p+1)\rs}{p[(p+1)\rs+q]}
	+ c_3 \cdot c_1^\frac{2}{p}
  \eas
  for all $t>0$.
\qed
We are now prepared for the announced testing procedure.
\begin{lem}\label{lem13}
  Suppose that (\ref{H1}) holds and that
  \be{13.1}
	\sup_{t>0} \io v^{\rs}(\cdot,t)<\infty
  \ee
  for some $\rs\ge 1$,
  and let $p\in (0,1)$ be such that $p<\frac{1}{\chi^2}$.
  Then with $\qpm$ taken from (\ref{1.1}), for all $q\in (\qm,\qp)$ fulfilling
  \be{13.11}
	q\le \frac{p(p+1)}{1-p}\cdot \rs
  \ee
  one can find $C>0$ such that
  \be{13.4}
	\io v^\frac{q}{p}(\cdot,t) \le C
	\qquad \mbox{for all } t>0.
  \ee
\end{lem}
\proof
  Since $\Omega$ is bounded, in view of (\ref{13.1}) it is sufficient to consider the case when $r:=\frac{q}{p}$
  satisfies $r>1$, and then testing the second equation in (\ref{0}) against $v^{r-1}$ shows that
  \bas
	\frac{1}{r} \frac{d}{dt} \io v^r
	+ (r-1) \io v^{r-2} v_x^2 + \io v^r
	= \io uv^r + \io B_2 v^{r-1}
	\qquad \mbox{for all } t>0.
  \eas
  Here, Young's inequality and the boundedness of $B_2$ show that there exists $c_1>0$ such that
  \bas
	\io B_2 v^{r-1} \le \frac{1}{2} \io v^r + c_1
	\qquad \mbox{for all } t>0,
  \eas
  so that $y(t):=\io v^r(\cdot,t)$, $t\ge 0$, satisfies
  \be{13.5}
	y'(t) + \frac{r}{2} y(t)
	\le h(t):=c_1 r + r\io u(\cdot,t)v(\cdot,t)
	\qquad \mbox{for all } t>0.
  \ee
  Now, thanks to our assumptions on $p$ and $q$, we may apply Lemma \ref{lem12} 
  to conclude from (\ref{13.1}) that there exists $c_2>0$ fulfilling
  \bas
	\int_t^{t+1} h(s)ds \le c_2
	\qquad \mbox{for all } t>0,
  \eas
  and therefore Lemma \ref{lem_ssw} ensures that (\ref{13.4}) is a consequence of (\ref{13.5}).
\qed
\subsection{Preparations for a recursive argument}
As Lemma \ref{lem13} suggests, our strategy toward improved estimates for $v$ will consist in a bootstrap-type
procedure, in a first step choosing $\rs:=1$ in Lemma \ref{lem13} and in each step seeking to maximize
the exponent $\frac{q}{p}$ appearing in (\ref{13.4}) according to our overall restrictions on $p$ and $q$
as well as (\ref{13.11}).
In order to create an appropriate framework for our iteration, let us introduce certain auxiliary functions,
and summarize some of their elementary properties, in the following lemma.
\begin{lem}\label{lem14}
  Let $\ps:=\min\{1,\frac{1}{\chi^2}\}$ as well as
  \be{14.1}
	\varphi_1(p):=\frac{p+1}{1-p},
	\quad
	\varphi_2(p):=\frac{1-p}{2p} \Big(1+\sqrt{1-p\chi^2}\Big)
	\quad \mbox{and} \quad
	\varphi_3(p):=\frac{1-p}{2p} \Big(1-\sqrt{1-p\chi^2}\Big)
   \ee
  for $p\in (0,\ps)$. Then 
  \be{14.2}
	\varphi_1'>0
	\quad \mbox{and} \quad
	\varphi_2'<0
	\qquad \mbox{on } (0,\ps),
  \ee
  and we have
  \be{14.3}
	\varphi_1(p) > \lim_{s\searrow 0} \varphi_1(s)=1
	\qquad \mbox{for all } p\in (0,\ps),
  \ee
  and
  \be{14.5}
	\varphi_2(p)\to + \infty
	\qquad \mbox{as } p\searrow 0,
  \ee
  as well as
  \be{14.4}
	\varphi_2(p) >\varphi_3(p)
	\qquad \mbox{for all } p\in (0,\ps).
  \ee
\end{lem}
\proof
  All statements can be verified by elementary computations.
\qed
Now the following observation explains the role of our smallness condition on $\chi$ from Theorem \ref{theo52}.
\begin{lem}\label{lem15}
  Suppose that $\chi<\frac{\sqrt{6\sqrt{3}+9}}{2}$.
  Then
  \be{15.1}
	\varphi_1(p_0)>\varphi_3(p_0)
  \ee
  is valid for the number
  \be{15.2}
	p_0:=\frac{2\sqrt{3}-3}{3} \in (0,1)
  \ee
  satisfying
  \be{15.3}
	p_0<\frac{1}{\chi^2}.
  \ee
\end{lem}
\proof
  We only need to observe that our assumption on $\chi$ warrants that
  \bas
	p_0\chi^2 < \frac{2\sqrt{3}-3}{3} \cdot \frac{6\sqrt{3}+9}{2} = \frac{3}{4},
  \eas
  which namely in particular yields (\ref{15.3}) and moreover implies that by (\ref{14.1}),
  \bas
	\frac{\varphi_3(p_0)}{\varphi_1(p_0)} -1
	&=& \frac{(1-p_0)^2}{2p_0(p_0+1)} \cdot \Big(1-\sqrt{1-p_0\chi^2}\Big) -1 \\
	&<& \frac{(1-p_0)^2}{2p_0(p_0+1)} \cdot \frac{1}{2} -1 \\[1mm]
	&=& 0,
  \eas
  as claimed.
\qed
Indeed, the latter property allows us to construct an increasing divergent sequence $(r_k)_{k\in\N}$ of exponents
to be used in Lemma \ref{lem13}.
\begin{lem}\label{lem151}
  Suppose that $\chi<\frac{\sqrt{6\sqrt{3}+9}}{2}$,
  and that $p_0$ is as in Lemma \ref{lem15}.
  Then for each $r\ge 1$, the set
  \be{151.1}
	S(r):=\Big\{ p\in (0,p_0) \ \Big| \ \varphi_2(p) \ge \varphi_1(p)\cdot r \Big\}
  \ee
  is not empty, and letting $r_0:=1$ as well as
  \be{151.2}
	p_k:=\sup S(r_r),
	\qquad k\in\N, 
  \ee
  and
  \be{151.3}
	r_k:=\varphi_1(p_k) \cdot r_{k-1},
	\qquad k\in\N,
  \ee
  recursively defines sequences $(p_k)_{k\in\N} \subset (0,p_0]$ and $(r_k)_{k\in\N} \subset (1,\infty)$ satisfying
  \be{151.4}
	p_k \le p_{k-1},
	\qquad \mbox{for all } k\in\N
  \ee
  and
  \be{151.5}
	r_k>r_{k-1},
	\qquad \mbox{for all } k\in\N
  \ee
  as well as
  \be{151.6}
	r_k\to\infty,
	\qquad \mbox{as } k\to\infty.
  \ee
  Moreover, writing
  \be{151.7}
	q_k:=p_k r_k,
	\qquad k\in\N,
  \ee
  we have
  \be{151.8}
	q^-(p_k) < q_k \le q^+(p_k)
	\qquad \mbox{for all } k\in\N
  \ee
  as well as
  \be{151.9}
	q_k \le \frac{p_k(p_k+1)}{1-p_k} \cdot r_{k-1}
	\qquad \mbox{for all } k\in\N.
  \ee
\end{lem}
\proof
  Observing that $\varphi_1$ and $\varphi_2$ are well-defined on $(0,p_0)$ due to the fact that 
  $p_0<\frac{1}{\chi^2} \le \ps$ by (\ref{15.3}), from (\ref{14.3}) and (\ref{14.5}) we see that
  \bas
	\frac{\varphi_2(p)}{\varphi_1(p)} \to + \infty
	\qquad \mbox{as } p\searrow 0,
  \eas
  implying that indeed $S(r)\ne\emptyset$ for all $r\ge 1$
  and that hence the definitions of $(p_k)_{k\in\N}$ and $(r_k)_{k\in\N}$ are meaningful.
  Moreover, from (\ref{151.2}) and (\ref{151.1}) it is evident that $p_k\in (0,p_0]$ for all $k\in\N$, whereas
  (\ref{151.3}) together with (\ref{14.3}) guarantees (\ref{151.5}) and that thus also the inclusion
  $(r_k)_{k\in\N} \subset (1,\infty)$ holds; as therefore $S(r_k) \subset S(r_{k-1})$ for all $k\in\N$, it is also
  clear that (\ref{151.4}) is valid.\\
  In order to verify (\ref{151.6}), assuming on the contrary that
  \be{151.99}
	r_k\to r_\infty
	\qquad \mbox{as } k\to\infty
  \ee
  with some $r_\infty \in (1,\infty)$, we would firstly obtain from (\ref{151.4}) that
  \be{151.10}
	p_k\searrow 0
	\qquad \mbox{as } k\to\infty,
  \ee
  for otherwise there would exist $p_\infty \in (0,p_0]$ such that $p_k\ge p_\infty$ for all $k\in\N$, which
  by (\ref{14.1}) would imply that $\varphi_1(p_k)\ge c_1:=\varphi_1(p_\infty)>1$ for all $k\in\N$ and that hence
  $r_k\ge c_1 r_{k-1}$ for all $k\in\N$ due to (\ref{151.3}), clearly contradicting the assumed boundedness property
  of $(r_k)_{k\in\N}$.
  In particular, (\ref{151.10}) entails the existence of $k_0\in\N$ such that
  \be{151.11}
	\varphi_2(p_k)=\varphi_1(p_k) \cdot r_{k-1}
	\qquad \mbox{for all } k\ge k_0,
  \ee
  because if this was false then for all $k\in\N$ we would have
  $\varphi_2(p)>\varphi_1(p)\cdot r_k$ for any $p\in (0,p_0)$ and thus $p_k=p_0$ for all $k\in\N$ by (\ref{151.2}).
  Now combining (\ref{151.11}) with (\ref{151.10}), however, again using (\ref{14.5}) we could infer that
  \bas
	\varphi_1(p_k) \cdot r_{k-1} = \varphi_2(p_k)\to + \infty
	\qquad \mbox{as } k\to\infty,
  \eas
  which is incompatible with the observation that
  \bas
	\varphi_1(p_k)\cdot r_{k-1} \to r_\infty<\infty
	\qquad \mbox{as } k\to\infty,
  \eas
  as asserted by (\ref{151.10}), (\ref{14.3}) and (\ref{151.99}).\abs
  To see that the numbers $q_k$ in (\ref{151.7}) have the claimed properties, we firstly use their definition along
  with those of $r_k$ and $\varphi_1$ to find that
  \bas
	q_k=p_k r_k = p_k \varphi_1(p_k) r_{k-1}
	= \frac{p_k(p_k+1)}{1-p_k} \cdot r_{k-1}
	\qquad \mbox{for all } k\in\N,
  \eas
  while from (\ref{151.2}) and (\ref{151.1}) it follows that $\varphi_1(p_k)\cdot r_{k-1} \le \varphi_2(p_k)$
  and thus
  \bas
	q_k= p_k \varphi_1(p_k) r_{k-1}
	\le p_k \varphi_2(p_k)
	= \frac{1-p_k}{2} \cdot \Big(1+\sqrt{1-p_k\chi^2}\Big)
	= q^+(p_k)
	\qquad \mbox{for all } k\in\N.
  \eas
  Finally, for the derivation of the left inequality in (\ref{151.8}) we make use of the property (\ref{15.1}) of $p_0$:
  Namely, if $k\in\N$ is such that $\varphi_2(p)\ge \varphi_1(p) \cdot r_k$ for all $p\in (0,p_0)$,
  then (\ref{151.2}) says that $p_k=p_0$ and therefore, by (\ref{151.7}), (\ref{151.3}), (\ref{151.5}), (\ref{15.1})
  and (\ref{14.1}),
  \bas
	q_k &=& p_k r_k = p_k \varphi_1(p_k) r_{k-1}
	= p_0 \varphi_1(p_0) r_{k-1} \ge p_0 \varphi_1(p_0) \\
	&>& p_0 \varphi_3(p_0) = p_k \varphi_3(p_k) 
	= \frac{1-p_k}{2}\Big( 1-\sqrt{1-p_k\chi^2}\Big)
	= q^-(p_k).
  \eas
  On the other hand, in the case when $k\in\N$ is such that
  $\inf_{p\in (0,p_0)} \big\{ \varphi_2(p)-\varphi_1(p)\cdot r_k\big\}$ is negative, (\ref{151.2}) implies that
  necessarily $\varphi_2(p_k)=\varphi_1(p_k)\cdot r_{k-1}$, so that
  \bas
	q_k=p_k \varphi_1(p_k) r_{k-1}
	= p_k \varphi_2(p_k)
	= \frac{1-p_k}{2}\Big(1+\sqrt{1-p_k\chi^2}\Big)
	> \frac{1-p_k}{2}\Big(1-\sqrt{1-p_k\chi^2}\Big),
  \eas
  because the restriction $p_k\le p_0$ together with (\ref{15.3}) ensures that $\sqrt{1-p_k \chi^2}$ must
  be positive.
\qed
\subsection{Boundedness of $v$ in $L^r(\Omega)$ for arbitrary $r<\infty$}
A straightforward induction on the basis of Lemma \ref{lem13} and Lemma \ref{lem151} leads to the following.
\begin{lem}\label{lem152}
  Let $\chi<\frac{\sqrt{6\sqrt{3}+9}}{2}$
  and suppose that \eqref{H1} holds, and let $(r_k)_{k\in\N_0} \subset (1,\infty)$ be as in Lemma \ref{lem151}.
  Then for all $k\in\N_0$ and any $r\in (1,r_k)\cup \{1\}$ there exists $C>0$ such that
  \be{152.1}
	\io v^r(\cdot,t) \le C
	\qquad \mbox{for all } t>0.
  \ee
\end{lem}
\proof
  Since for $k=0$ this has been asserted by Lemma \ref{lem02}, in view of an inductive argument we only need to 
  make sure that if for some $k\in\N$ we have
  \be{152.2}
	\sup_{t>0} \io v^r(\cdot,t) <\infty
	\qquad \mbox{for all } r\in (1,r_{k-1}) \cup \{1\},
  \ee
  then 
  \be{152.3}
	\sup_{t>0} \io v^r(\cdot,t) <\infty
	\qquad \mbox{for all } r\in (1,r_k).
  \ee
  In verifying this, 
  by boundedness of $\Omega$ we may concentrate on values of $r\in (1,r_k)$ which are sufficiently close
  to $r_k$ such that with $p_k$ as in Lemma \ref{lem151} and $q^-(p_k)$ taken from (\ref{1.1})
  we have
  \be{152.4}
	r>\frac{q^-(p_k)}{p_k},
  \ee
  which is possible since from (\ref{151.7}) and (\ref{151.8}) we know that
  \bas
	p_k r \to q_k= p_k r_k>q^-(p_k)
	\qquad \mbox{as } r\to r_k.
  \eas
  We now let
  \be{152.44}
	q:=p_k r
  \ee
  and
  \be{152.5}
	\rs:=\max \Big\{ 1 \, , \, \frac{(1-p_k)q}{p_k(p_k+1)} \Big\}
  \ee
  and observe that then
  \be{152.6}
	q>q^-(p_k)
  \ee
  by (\ref{152.4}) and
  \be{152.7}
	q< p_k r_k \le q^+(p_k)
  \ee
  by (\ref{151.7}) and (\ref{151.8}), whereas (\ref{152.5}) ensures that
  \be{152.8}
	q \le \frac{p_k(p_k+1)}{1-p_k} \cdot \rs.
  \ee
  From (\ref{152.5}) it moreover follows that if $\rs>1$ then since $r<r_k$ implies that $q<q_k$, we have
  \bas
	\rs=\frac{(1-p_k)q}{p_k(p_k+1)}
	<\frac{(1-p_k)q_k}{p_k(p_k+1)} \le r_{k-1}
  \eas
  according to (\ref{151.9}).
  As thus (\ref{152.2}) warrants that
  \bas
	\sup_{t>0} \io v^{\rs}(\cdot,t)<\infty,
  \eas
  in view of (\ref{152.6}), (\ref{152.7}) and (\ref{152.8}) we may apply Lemma \ref{lem13} to find $c_1>0$ such that
  \bas
	\io v^\frac{q}{p_k}(\cdot,t) \le c_1
	\qquad \mbox{for all } t>0,
  \eas
  which thanks to (\ref{152.44}) yields (\ref{152.3}), because $r$ was an arbitrary number in the range described
  in (\ref{152.3}) and (\ref{152.4}).
\qed
In particular, $v$ remains bounded in $L^r(\Omega)$ for arbitrarily large finite $r$:
\begin{cor}\label{cor153}
  Let $\chi<\frac{\sqrt{6\sqrt{3}+9}}{2}$, and assume \eqref{H1}.
  Then for all $r\ge 1$ there exists $C>0$ such that
  \bas
	\io v^r(\cdot,t) \le C
	\qquad \mbox{for all } t>0.
  \eas
\end{cor}
\proof
  Since Lemma \ref{lem151} asserts that the sequence $(r_k)_{k\in\N}$ introduced there has the property that
  $r_k\to\infty$ as $k\to\infty$, this is an immediate consequence of Lemma \ref{lem152}.
\qed
\subsection{H\"older regularity of $v$}
Once more relying on the first-order estimate provided by Lemma \ref{lem10} and the basic property 
$\int_0^\infty \io uv<0$ asserted by Lemma \ref{lem01}, from Corollary \ref{cor153} we can now derive
boundedness, and even a certain temporal decay,
of the forcing term $uv$ from the second equation in (\ref{0}) with respect to some superquadratic 
space-time Lebesgue norm.
\begin{lem}\label{lem17}
  Let $\chi<\frac{\sqrt{6\sqrt{3}+9}}{2}$, and assume \eqref{H1}.
  Then for all $p\in (0,\frac{1}{3})$ fulfilling $p<\frac{1}{\chi^2}$ we have
  \be{17.1}
	\int_t^{t+1} \io (uv)^{p+2} \to 0
	\qquad \mbox{as } t\to\infty.
  \ee
\end{lem}
\proof
  We first note that taking $\xi\to\infty$ in Lemma \ref{lem01} shows that our hypothesis (\ref{H1}) 
  warrants that $\int_0^\infty \io uv<\infty$ and hence
  \bas
	\int_t^{t+1} \io uv \to 0
	\qquad \mbox{as } t\to\infty.
  \eas
  In view of an interpolation argument, it is therefore sufficient to make sure that for all $\tp \in (0,\frac{1}{3})$
  satisfying $\tp<\frac{1}{\chi^2}$ we can find $c_0>0$ such that
  \be{17.2}
	\int_t^{t+1} \io (uv)^{\tp+2} \le c_1
	\qquad \mbox{for all } t>0.
  \ee
  For this purpose, given any such $\tp$ we can fix $p\in (\tp,\frac{1}{3})$ such that still $p<\frac{1}{\chi^2}$,
  and then observe that
  \bas
	\frac{3p-1}{1-p} < 0 < \sqrt{1-p\chi^2}.
  \eas
  This ensures that with the numbers $\qpm$ from (\ref{1.1}) we have $\qp>p$, whence it is possible to pick
  $q\in (\qm,\qp)$ such that $q>p$.
  Writing $r:=\tp+2$, by means of the H\"older inequality we can thus estimate
  \bas
	\io (uv)^r
	&=& \io \Big( u^\frac{p}{2} v^\frac{q}{2}\Big)^\frac{2r}{q} \cdot u^\frac{(q-p)r}{q} \\
	&\le& \bigg\{ \io \Big(u^\frac{p}{2} v^\frac{q}{2}\Big)^\frac{2r}{q-(q-p)r} \bigg\}^\frac{q-(q-p)r}{q}
	\cdot \bigg\{ \io u \bigg\}^\frac{(q-p)r}{q} \\
	&\le& c_2 \bigg\{ \io \Big(u^\frac{p}{2} v^\frac{q}{2}\Big)^\frac{2r}{q-(q-p)r} \bigg\}^\frac{q-(q-p)r}{q} \\
	&=& c_2 \Big\|u^\frac{p}{2} v^\frac{q}{2}\Big\|_{L^\frac{2r}{q-(q-p)r}(\Omega)}^\frac{2r}{q}
	\qquad \mbox{for all } t>0
  \eas
  with $c_2:=\sup_{t>0} \|u(\cdot,t)\|_{L^1(\Omega)}^\frac{(q-p)r}{q}$ being finite according to 
  Lemma \ref{lem01} and our assumption that (\ref{H1}) be valid.\\
  Consequently, using the Gagliardo-Nirenberg inequality we see that with some $c_3>0$ we have
  \bea{17.3}
	\int_t^{t+1} \io (uv)^r
	&\le& c_3 \int_t^{t+1} \Big\| \Big(u^\frac{p}{2}(\cdot,s) v^\frac{q}{2}(\cdot,s)\Big)_x \Big\|_{L^2(\Omega)}^2
	\Big\| u^\frac{p}{2}(\cdot,s) v^\frac{q}{2}(\cdot,s)\Big\|_{L^\frac{2}{p+\eps q}(\Omega)}^\frac{2(r-q)}{q} ds 
		\nn\\
	& & + c_3 \int_t^{t+1} 
	\Big\| u^\frac{p}{2}(\cdot,s) v^\frac{q}{2}(\cdot,s)\Big\|_{L^\frac{2}{p+\eps q}(\Omega)}^\frac{2r}{q} ds
	\qquad \mbox{for all } t>0,
  \eea
  where we have abbreviated
  \bas
	\eps:=\frac{p+2-r}{r-q}.
  \eas
  Now since $\eps$ is positive because $r<p+2$ and $r>2>1>\qp>q$, and since thus $\frac{2}{p+\eps q}<\frac{2}{p}$,
  an application of Lemma \ref{lem11} readily yields $c_4>0$ such that
  \bas
	\Big\| u^\frac{p}{2}(\cdot,s) v^\frac{q}{2}(\cdot,s)\Big\|_{L^\frac{2}{p+\eps q}(\Omega)}
	\le c_4
	\qquad \mbox{for all } s>0,
  \eas
  whereas the inequalities $p<\min \{1,\frac{1}{\chi^2}\}$ and $\qm<q<\qp$ ensure that due to Lemma \ref{lem10}
  we can find $c_5>0$ such that
  \bas
	\int_t^{t+1} \Big\| \Big(u^\frac{p}{2}(\cdot,s) v^\frac{q}{2}(\cdot,s)\Big)_x \Big\|_{L^2(\Omega)}^2 ds
	\le c_5
	\qquad \mbox{for all } t>0.
  \eas
  Therefore, (\ref{17.3}) implies that
  \bas
	\int_t^{t+1} \io (uv)^r \le c_3 c_4^\frac{2(r-q)}{q} c_5 + c_3 c_4^\frac{2r}{q}
	\qquad \mbox{for all } t>0
  \eas
  and hence proves (\ref{17.2}) due to our definition of $r$.
\qed
Thanks to the fact that the integrability exponent appearing therein is large than $2$, 
the boundedness property implied by the decay statement in Lemma \ref{lem17} allows us to derive 
boundedness of $v$ even in a space compactly embedded into $L^\infty(\Omega)$.
\begin{lem}\label{lem18}
  Let $\chi<\frac{\sqrt{6\sqrt{3}+9}}{2}$ and assume \eqref{H1}.
  Then there exist $\gamma\in (0,1)$ and $C>0$ such that
  \be{18.1}
	\|v(\cdot,t)\|_{C^\gamma(\bom)} \le C
	\qquad \mbox{for all $t>1$.}
  \ee
\end{lem}
\proof
  We fix $\beta \in (\frac{1}{4},\frac{1}{2})$ and any $\gamma\in (0,2\beta-\frac{1}{2})$ and then once more refer to
  known embedding results (\cite{henry}) to recall that the sectorial realization $A$ of $-(\cdot)_{xx}+1$ 
  under homogeneous 
  Neumann boundary conditions in $L^2(\Omega)$ has the property that its fractional power $A^\beta$ satisfies
  $D(A^\beta) \hra C^\gamma(\bom)$. Therefore, writing
  \bas
	v(\cdot,t)=e^{-A} v(\cdot,t-1) + \int_{t-1}^t e^{-(t-s)A} h(\cdot,s) ds
	\qquad \mbox{for } t>1
  \eas
  with
  \bas
	h(\cdot,t):=u(\cdot,t)v(\cdot,t)+B_2(\cdot,t),
	\qquad t>0,
  \eas
  we can estimate
  \bas
	\|v(\cdot,t)\|_{C^\gamma(\bom)} 
	\le c_1 \Big\|A^\beta e^{-{\color{red}{t}}A} v(\cdot,t-1)\Big\|_{L^2(\Omega)}
	+ c_1 \int_{t-1}^t \Big\| A^\beta e^{-(t-s)A} h(\cdot,s)\Big\|_{L^2(\Omega)} ds,
	\qquad \mbox{for all } t>1,
  \eas
  with some $c_1>0$.
  As well-known regularization features of $(e^{-tA})_{t\ge 0}$ (\cite{friedman}) warrant the existence of $c_2>0$ 
  fulfilling
  \bas
	\Big\| A^\beta e^{-tA} \varphi\Big\|_{L^2(\Omega)}
	\le c_2 t^{-\beta} \|\varphi\|_{L^2(\Omega)}
	\qquad \mbox{for all $\varphi\in C^0(\bom)$ and any } t>0,
  \eas
  by using the Cauchy-Schwarz inequality we infer that
  \bea{18.2}
	\|v(\cdot,t)\|_{C^\gamma(\bom)}
	&\le& c_1 c_2\|v(\cdot,t-1)\|_{L^2(\Omega)}
	+ c_1 c_2 \int_{t-1}^t (t-s)^{-\beta} \|h(\cdot,s)\|_{L^2(\Omega)} ds \nn\\
	&\le& c_1 c_2\|v(\cdot,t-1)\|_{L^2(\Omega)}
	+ c_1 c_2 \bigg\{ \int_{t-1}^t (t-s)^{-2\beta} ds \bigg\}^\frac{1}{2}
	\bigg\{ \int_{t-1}^t \|h(\cdot,s)\|_{L^2(\Omega)}^2 ds \bigg\}^\frac{1}{2}
  \eea
  for all $t>1$, where we note that
  \bas
	\int_{t-1}^t (t-s)^{-2\beta} ds = \frac{1}{1-2\beta}
	\qquad \mbox{for all } t>1
  \eas
  thanks to our restriction $\beta<\frac{1}{2}$.
  Since Corollary \ref{cor153} provides $c_3>0$ such that
  \bas
	\|v(\cdot,t-1)\|_{L^2(\Omega)} \le c_3
	\qquad \mbox{for all } t>1,
  \eas
  and since Lemma \ref{lem17} along with the boundedness of $B_2$ in $\Omega\times (0,\infty)$ implies that
  \bas
	\int_{t-1}^t \|h(\cdot,s)\|_{L^2(\Omega)}^2 ds \le c_4
	\qquad \mbox{for all } t>1
  \eas
  with some $c_4>0$, the inequality in (\ref{18.1}) is thus a consequence of (\ref{18.2}).
\qed
\subsection{Stabilization of $v$ under the hypotheses (\ref{H1}) and (\ref{H3}). Proof of Theorem \ref{theo52}}
\begin{figure}[H]
  \center
      \subfloat[$\chi=20$]{\label{fig:no_env}\includegraphics[width=0.35\textwidth]{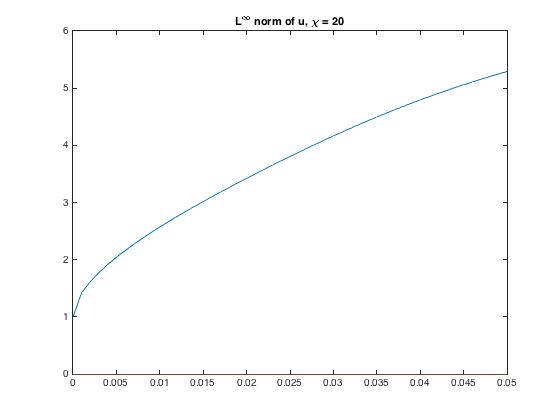}}\;
     \subfloat[$\chi=50$]{\label{fig:no_env}\includegraphics[width=0.35\textwidth]{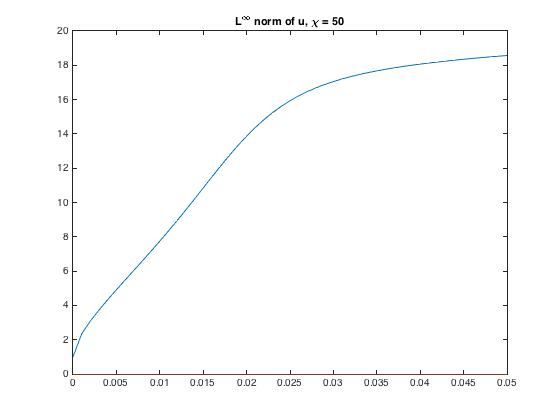}}\\
  \subfloat[$\chi=100$]{\label{fig:no_env}\includegraphics[width=0.35\textwidth]{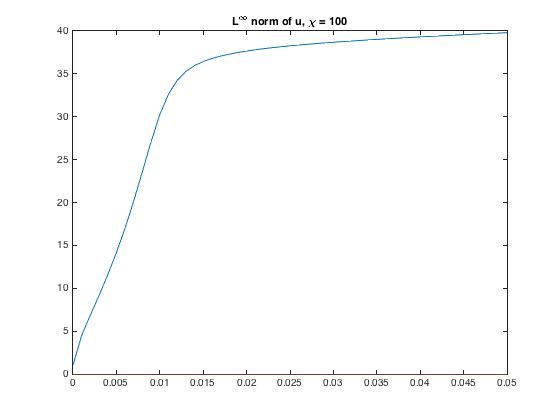}}\;
  \subfloat[$\chi=150$]{\label{fig:no_env}\includegraphics[width=0.35\textwidth]{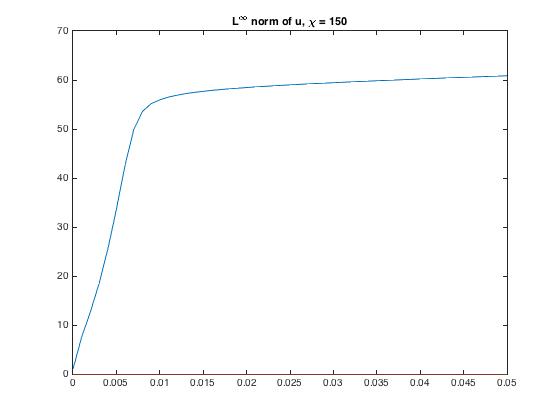}}\\
    \subfloat[$\chi=500$]{\label{fig:no_env}\includegraphics[width=0.35\textwidth]{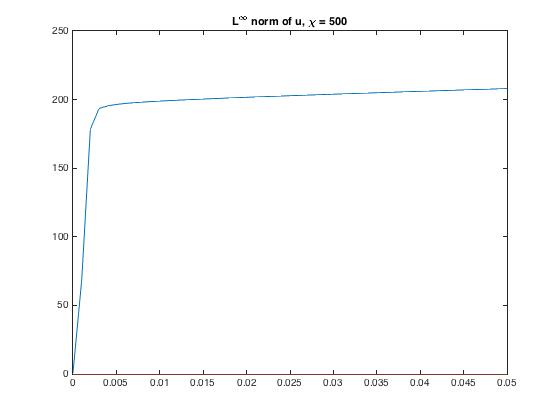}}\;
  \subfloat[$\chi=1000$]{\label{fig:no_env}\includegraphics[width=0.35\textwidth]{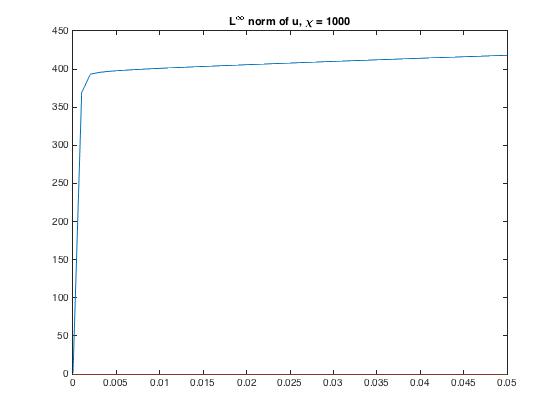}}
   \caption{The evolution of the maximum concentration of criminal $\norm{u(\cdot,t)}_{\infty}$ 
   at a short time scale $t\in [0,.05]$ with initial condition given by $(u(x,0),v(x,0))=(e^{-x},e^{-x})$ and $B_1 = B_2 =1$.  }   \label{fig:short}
\end{figure}

As a final preparation for the proof of Theorem \ref{theo52}, let us now make use of the $L^2$ decay property
of $uv$ entailed by Lemma \ref{lem17} in order to assert that under the additional assumption (\ref{H3}),
$v$ indeed stabilizes toward the desired limit, at least with respect to the topology in $L^2(\Omega)$.
\begin{lem}\label{lem19}
  Let $\chi<\frac{\sqrt{6\sqrt{3}+9}}{2}$, and assume that \eqref{H1} and \eqref{H3} hold with some
  $B_{2,\infty} \in L^2(\Omega)$.
  Then 
  \be{19.2}
	v(\cdot,t) \to v_\infty
	\quad \mbox{in } L^2(\Omega)
	\qquad \mbox{as } t\to\infty,
  \ee
  where $v_\infty$ denotes the solution of \eqref{vinfty}.
\end{lem}
\proof
  Using (\ref{0}) and (\ref{vinfty}) we compute
  \bas
	\frac{1}{2} \frac{d}{dt} \io (v-v_\infty)^2
	&=& \io (v-v_\infty)\cdot (v_{xx}+uv-v+B_2) \\
	&=& \io (v-v_\infty)\cdot \Big\{ (v-v_\infty)_{xx} - (v-v_\infty) + uv + (B_2-B_{2,\infty}) \Big\} \\
	&=& - \io (v-v_\infty)_x^2 - \io (v-v_\infty)^2
	+ \io (v-v_\infty) \cdot \Big\{ uv+ (B_2-B_{2,\infty}) \Big\}
	\qquad \mbox{for all } t>0,
  \eas
  where the first summand on the right is nonpositive, and where the rightmost integral can be estimated by Young's
  inequality according to
  \bas
	\io (v-v_\infty) \cdot \Big\{ uv+ (B_2-B_{2,\infty}) \Big\}
	&\le& \frac{1}{2} \io (v-v_\infty)^2
	+ \frac{1}{2} \io \Big\{ uv+(B_2-B_{2,\infty}) \Big\}^2 \\
	&\le& \frac{1}{2} \io (v-v_\infty)^2
	+ \io (uv)^2 + \io (B_2-B_{2,\infty})^2
	\qquad \mbox{for all } t>0.
  \eas
  Therefore, $y(t):=\io (v(\cdot,t)-v_\infty)^2$ and
  $h(t):=2\io (u(\cdot,t)v(\cdot,t))^2 + 2\io (B_2(\cdot,t)-B_{2,\infty})^2$, \ $t\ge 0$, satisfy
  \bas
	y'(t)+y(t) \le h(t)
	\qquad \mbox{for all } t>0,
  \eas
  so that since Lemma \ref{lem17} entails that
  \bas
	\int_t^{t+1} \io (uv)^2 \to 0
	\qquad \mbox{as } t\to\infty
  \eas
  and that thus
  \bas
	\int_t^{t+1} h(s)ds \to 0
	\qquad \mbox{as } t\to\infty
  \eas
  thanks to (\ref{H3}), the claimed property (\ref{19.2}) results from Lemma \ref{lem_espejo_win}.
\qed
Collecting all the above, we can easily derive our main result on asymptotic behavior under the assumptions
that (\ref{H1}) and possibly also (\ref{H3}) hold.\abs
\proofc of Theorem \ref{theo52}.\quad
  Supposing that $\chi < \frac{\sqrt{6\sqrt{3}+9}}{2}$ and that (\ref{H1}) be valid,
  we obtain the boundedness property (\ref{52.2}) of $v$ in $\Omega\times (0,\infty)$ as a consequence 
  of Lemma \ref{lem18} and Lemma \ref{lem_loc}.
  If moreover (\ref{H3}) is fulfilled with some $B_{2,\infty}\in L^2(\Omega)$, then from 
  Lemma \ref{lem19} we know that $v(\cdot,t)\to v_\infty$ in $L^2(\Omega)$ as $t\to\infty$.
  Since Lemma \ref{lem18} actually even warrants precompactness of $(v(\cdot,t))_{t>1}$ in $L^\infty(\Omega)$
  by means of the Arzel\`a-Ascoli theorem, this already implies the uniform convergence property 
  claimed in (\ref{52.3}).
\qed
\mysection{Numerical results}
In this section we explore the growth of solutions to \eqref{0} as $\chi$ increases on 
small time scales.  The effect of large chemotaxis sensitivities on the growth of the solutions
has been observed in Keller-Segel-type systems.  From numerical simulations we observe that
the $L^\infty$ norm of the criminal density increases sharply with $\chi$ in short-time scales before relaxing
to the steady-state solution.  Indeed, the solution quickly relaxes to a steady-state solution
once the dissipation is able to dominate.  For all numerical experiments we consider
initial data $u(x,0)= e^{-x}$ and $v(x,0)= e^{-x},$ $B_1 = B_2 = 1$ and vary the parameter $\chi$.  
All numerical computations where made using Matlab's {\it pdepe} function.  In
Figure \ref{fig:short} we observe the rapid growth on the short time scale ($t \in [0,.05]$ with time step $\delta t = .001$).
This figure illustrates the fact that the criminal density reaches a higher value as $\chi$ increases.  

\begin{figure}[H]
  \center
      \subfloat[$\chi=20$]{\label{fig:no_env}\includegraphics[width=0.35\textwidth]{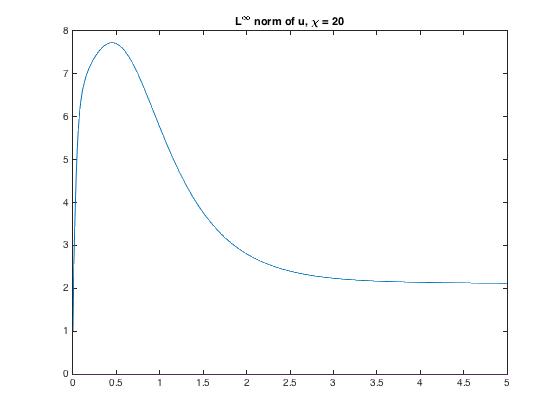}}\;
     \subfloat[$\chi=50$]{\label{fig:no_env}\includegraphics[width=0.35\textwidth]{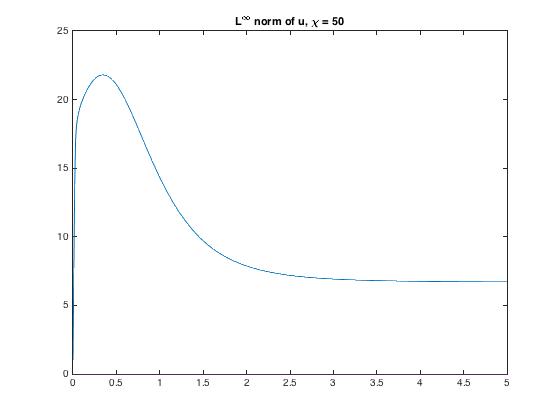}}\\
  \subfloat[$\chi=100$]{\label{fig:no_env}\includegraphics[width=0.35\textwidth]{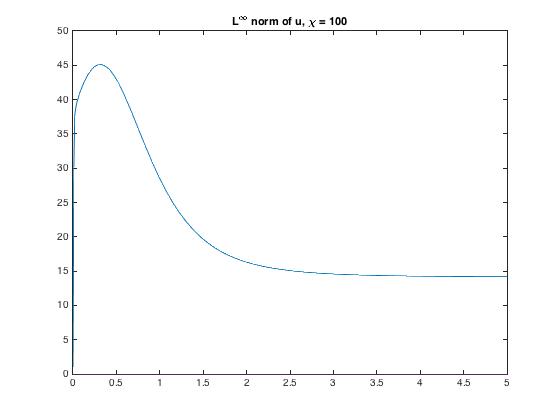}}\;
  \subfloat[$\chi=150$]{\label{fig:no_env}\includegraphics[width=0.35\textwidth]{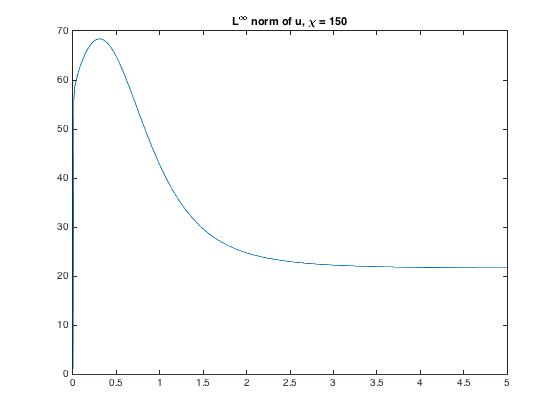}}\\
    \subfloat[$\chi=500$]{\label{fig:no_env}\includegraphics[width=0.35\textwidth]{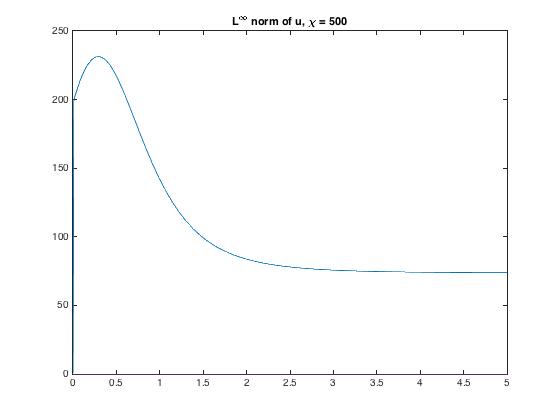}}\;
  \subfloat[$\chi=1000$]{\label{fig:no_env}\includegraphics[width=0.35\textwidth]{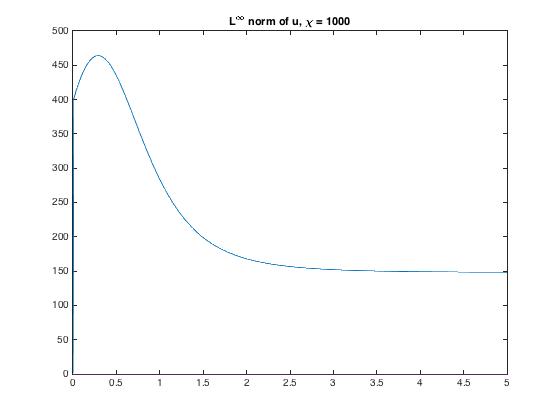}}
   \caption{The evolution of the maximum concentration of criminal $\norm{u(\cdot,t)}_{\infty}$ 
   at a longer time scale ($t\in [0,5]$) with initial condition given by $(u(x,0),v(x,0))=(e^{-x},e^{-x})$ and $B_1 = B_2 =1$.  }   \label{fig:long}
\end{figure}

\begin{figure}[H]
  \center
      \subfloat[$\chi=12$]{\label{fig:no_env}\includegraphics[width=0.35\textwidth]{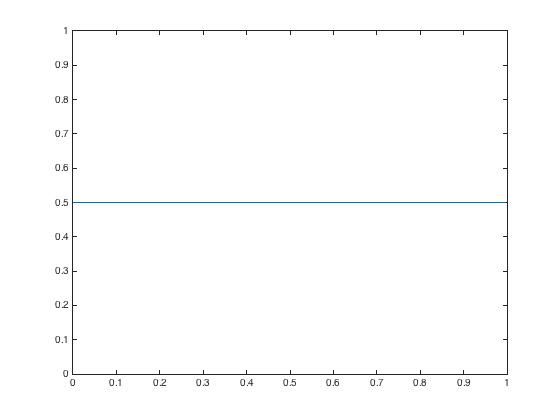}}\;
    \subfloat[$\chi=13$]{\label{fig:no_env}\includegraphics[width=0.35\textwidth]{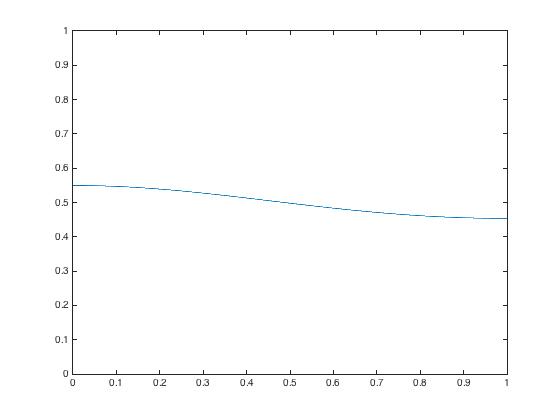}}\\
  \subfloat[$\chi=20$]{\label{fig:no_env}\includegraphics[width=0.35\textwidth]{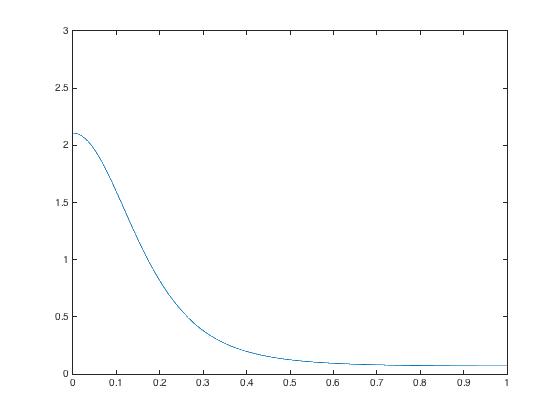}}\;
  \subfloat[$\chi=50$]{\label{fig:no_env}\includegraphics[width=0.35\textwidth]{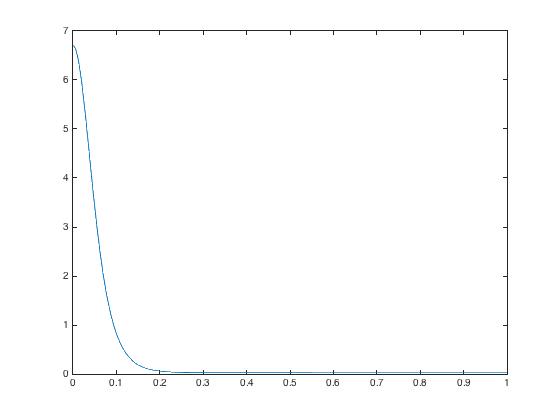}}\\
   \caption{Criminal density $u(x,t)$ at $t =20$ for various values of $\chi$.  }   \label{fig:soln}
\end{figure}

On the other hand, at longer time scales (although not so long $t\in [0,5]$ with $\delta t = .05$) the dissipation dominates and in all cases we see eventual decay.  This is
illustrated in Figure \ref{fig:long}, where we can see that by time $t=5$ the maximum density of criminals
has reached a steady state.  
Another interesting thing to note that is that the steady-state of the maximum density of criminals increases with $\chi$.  
Thus, we do not see a relaxation to the constant steady states, which in this case are $u\equiv \frac{1}{2}$ and $v \equiv 2$.  
In fact, relaxation to the homogeneous steady-states occurs with $\chi$ small.  However, as $\chi$ increases we
observe to a non-constant hump-solution with the maximum at the origin.  Figure \ref{fig:soln} 
\mysection{Appendix: Two ODE lemmata}
Let us separately formulate two auxiliary statements on boundedness and decay in linear
ODIs with inhomogeneities enjoying certain averaged boundedness and decay properties.
\begin{lem}\label{lem_ssw}
  Let $T\in (0,\infty]$ and $\tau\in (0,T)$, and let $y\in C^1([0,T))$ and $h\in L^1_{loc}([0,\infty))$ 
  be nonnegative and such that with some $a>0$ and $b>0$ we have 
  \bas
	y'(t) + ay(t) \le h(t)
	\qquad \mbox{for all } t\in (0,T),
  \eas
  as well as
  \bas
	\frac{1}{\tau} \int_t^{t+\tau} h(s)ds \le b
	\qquad \mbox{for all } t\in (0,T).
  \eas
  Then
  \bas
	y(t) \le y(0) + \frac{b\tau}{1-e^{-a\tau}}
	\qquad \mbox{for all } t\in [0,T).
  \eas
\end{lem}
\proof 
  This can be found e.g.~in \cite[Lemma 3.4]{win_ks_nasto_exist}.
\qed
\begin{lem}\label{lem_espejo_win}
  Let $y\in C^1([0,\infty))$ and $h\in L^1_{loc}([0,\infty))$ be nonnegative functions satisfying
  \bas
	y'(t) + a y(t) \le h(t)
	\qquad \mbox{for all } t>0
  \eas
  with some $a>0$. Then if
  \bas
	\int_t^{t+1} h(s) ds \to 0
	\qquad \mbox{as } t\to\infty,
  \eas
  we have
  \bas
	y(t) \to 0
	\qquad \mbox{as } t\to\infty.
  \eas
\end{lem}
\proof
  An elementary derivation of this has been given in \cite[Lemma 4.6]{espejo_win1}, for instance.
\qed
\vspace*{5mm}
{\bf Acknowledgement.} \quad
  The second author acknowledges support of the Deutsche  
Forschungsgemeinschaft in the context of the project
{\em Emergence of structures and advantages in cross-diffusion  
systems} (No.-- 411007140, GZ: WI 3707/5-1).  The first author was supported by the
  NSF Grant DMS-1516778. 
\bibliographystyle{plain}
\bibliography{crime}

\begin{thebibliography}{10}

\bibitem{amann}
H.~Amann.
\newblock {Dynamic theory of quasilinear parabolic systems III. Global
  existence.}
\newblock {\em Math. Z.}, 202:219--250, 1989.

\bibitem{Barbaro2013}
A.~B.~T. Barbaro, L.~Chayes, and M.~R. D'Orsogna.
\newblock {Territorial developments based on graffiti: A statistical mechanics
  approach}.
\newblock {\em Physica A: Statistical Mechanics and its Applications},
  392(1):252--270, 2013.

\bibitem{BBTW}
N.~Bellomo, A~Bellouquid, Y.~Tao, and M.~Winkler.
\newblock Toward a mathematical theory of keller-segel models of pattern
  formation in biological tissues.
\newblock {\em Math. Mod.~Meth. Appl.~Sci.}, 25:1663--1763, 2015.

\bibitem{Berestycki2010}
H.~Berestycki and J-P. Nadal.
\newblock {Self-organised critical hot spots of criminal activity}.
\newblock {\em European J. Appl. Math.}, 21(4-5):371--399, 2010.

\bibitem{Berestycki2013c}
H.~Berestycki, N.~Rodr{\'{i}}guez, and L.~Ryzhik.
\newblock {Traveling wave solutions in a reaction-diffusion model for criminal
  activity}.
\newblock {\em Multiscale Modeling {\&} Simulation}, 11(4):1097--1126, 2013.

\bibitem{Berestycki2014a}
H.~Berestycki, J.~Wei, and M.~Winter.
\newblock {Existence of symmetric and asymmetric spikes for a crime hotspot
  model}.
\newblock {\em SIAM J. Math. Anal.}, 46(1):691--719, 2014.

\bibitem{biler}
P.~Biler.
\newblock {Global solutions to some parabolic-elliptic systems of chemotaxis}.
\newblock {\em Adv. Math. Sci. Appl.}, 9(1):347--359, 1999.

\bibitem{Cantrell2012}
R.~S. Cantrell, C.~Cosner, and R.~Man{\'{a}}sevich.
\newblock {Global bifurcation of solutions for crime modeling equations}.
\newblock {\em SIAM J. Appl. Math.}, 44(3):1340--1358, 2012.

\bibitem{Chaturapruek2013}
S.~Chaturapruek, J.~Breslau, D.~Yazidi, T.~Kolokolnikiv, and S.~McCalla.
\newblock {Crime modeling with Levy flights}.
\newblock {\em SIAM J. Appl. Math.}, 73(4):1703--1720, 2013.

\bibitem{Cohen1979}
L.~E. Cohen and M.~Felson.
\newblock {Social change and crime rate trends: a routine activity approach}.
\newblock {\em American Sociological Review}, 44(4):588--608, 1979.

\bibitem{DOrsogna2015}
M.~R. D'Orsogna and M.~Perc.
\newblock {Statistical physics of crime: A review}.
\newblock {\em Physics of Life Reviews}, 12:1--21, 2015.

\bibitem{espejo_win1}
E.~Espejo and M.~Winkler.
\newblock {Global classical solvability and stabilization in a two-dimensional
  chemotaxis-Navier-Stokes system modeling coral fertilization.}
\newblock {\em Preprint}.

\bibitem{Felson1987}
M.~Felson.
\newblock {Routine activities and crime prevention in the developing
  metropolis}.
\newblock {\em Criminology}, 25:911--932, 1987.

\bibitem{friedman}
A.~Friedman.
\newblock {\em {Partial Differential Equations}}.
\newblock Holt, Rinehart $\backslash${\&} Winston, New York, 1969.

\bibitem{Gu2016}
Y.~Gu, Q.~Wang, and Y.~Guangzeng.
\newblock {Stationary patterns and their selection mechanism of Urban crime
  models with heterogeneous near--repeat victimization effec}.
\newblock {\em Arxiv preprint:1409.0835v2}, pages 1--30, 2016.

\bibitem{henry}
D.~Henry.
\newblock {\em {Geometric Theory of Semilinear Parabolic Equations}}.
\newblock Springer Verlag, Berlin, 1981.

\bibitem{herrero_velazquez}
M.A. Herrero and J.J.L. Vel\'azquez.
\newblock A blow-up mechanism for a chemotaxis model.
\newblock {\em Ann.~Scuola Normale Superiore Pisa}, 24:633--683, 1997.

\bibitem{hpw}
T.~Hillen, K.~J. Painter, and M.~Winkler.
\newblock {Convergence of a Cancer Invasion Model To a Logistic Chemotaxis
  Model}.
\newblock {\em Mathematical Models and Methods in Applied Sciences},
  23(01):165--198, 2013.

\bibitem{horstmann_win}
D.~Horstmann and M.~Winkler.
\newblock {Boundedness vs. blow-up in a chemotaxis system}.
\newblock {\em Journal of Differential Equations}, 215(1):52--107, 2005.

\bibitem{Johnson1997}
S.~D. Johnson, K.~Bowers, and A.~Hirschfield.
\newblock {New insights into the spatial and temporal distribution of repeat
  victimisation}.
\newblock {\em British Journal of Criminology}, 37(2):224--241, 1997.

\bibitem{Jones2010}
P.~A. Jones, P.~J. Brantingham, and L.~R. Chayes.
\newblock {Statistical models of criminal behavior: the effects of law
  enforcement actions}.
\newblock {\em Mathematical Models {\&} Methods in Applied Sciences},
  20(1):1397--1423, 2010.

\bibitem{Kelling1982}
G.~L. Kelling and J.~Q. Wilson.
\newblock {Broken windows}, 1982.

\bibitem{Kolokolnikiv2014}
T.~Kolokolnikiv, M.~J. Ward, and J.~Wei.
\newblock {The stability of hotspot patterns for reaction-diffusion models of
  urban crime}.
\newblock {\em DCDS-B}, 19:1373--1410., 2014.

\bibitem{lankeit_MMAS}
J.~Lankeit.
\newblock A new approach toward boundedness in a two-dimensional parabolic
  chemotaxis system with singular sensitivity.
\newblock {\em Math.~Meth.~Appl.~Sci.}, 39:394 -- 404, 2016.

\bibitem{lw2}
J.~Lankeit and M.~Winkler.
\newblock A generalized solution concept for the keller--segel system with
  logarithmic sensitivity: global solvability for large nonradial data.
\newblock 4(24):24--49, 2017.

\bibitem{Manasevich2012}
R.~Man{\'{a}}sevich, Q.~H. Phan, and P.~Souplet.
\newblock {Global existence of solutions for a chemotaxis-type system arising
  in crime modelling}.
\newblock {\em European J. Appl. Math.}, 24(02):273--296, 2012.

\bibitem{McMillon2014}
D.~McMillon, C.~P. Simon, and J.~Morenoff.
\newblock {Modeling the underlying dynamics of the spread of crime}.
\newblock {\em PLoS ONE}, 9(4), 2014.

\bibitem{Nuno2011}
J.~C. Nu{\~{n}}o, M.~A. Herrero, and M.~Primicerio.
\newblock {A mathematical model of a criminal-prone society}.
\newblock {\em DCDS-S}, 4(1):193--207, 2011.

\bibitem{osaki_yagi}
K.~Osaki and A.~Yagi.
\newblock Finite dimensional attractor for one-dimensional keller-segel
  equations.
\newblock {\em Funkcialaj Ekvacioj}, 22:441--469, 2001.

\bibitem{Ricketson2010}
L.~Ricketson.
\newblock {A continuum model of residential burglary incorporating law
  enforcement}.
\newblock {\em Preprint}, pages 1--7, 2010.

\bibitem{Rodriguez2013}
N.~Rodr{\'{i}}guez.
\newblock {On the global well-posedness theory for a class of PDE models for
  criminal activity}.
\newblock {\em Physica D: Nonlinear Phenomena}, 260(3):191--200, oct 2013.

\bibitem{Rodriguez2010}
N.~Rodr{\'{i}}guez and A.L. Bertozzi.
\newblock {Local existence and uniqueness of solutions to a PDE model for
  criminal behavior}.
\newblock {\em Mathematical Models and Methods in Applied Sciences},
  20(supp01):1425--1457, sep 2010.

\bibitem{Rodriguez2014a}
N.~Rodr{\'{i}}guez and L.~Ryzhik.
\newblock {Exploring the effects of social preference, economic disparity, and
  heterogeneous environments on segregation}.
\newblock {\em Communications in Mathematical Sciences}, 14(2):363--387, 2016.

\bibitem{Short2010}
M.~B. Short, A.L. Bertozzi, and P.~J. Brantingham.
\newblock {Nonlinear patterns in urban crime: hotspots, bifurcations, and
  suppression}.
\newblock {\em SIAM Journal on Applied Dynamical Systems}, 9(2):462--483, 2010.

\bibitem{Short2009}
M.~B. Short, M.~R. D'Orsogna, P.~J. Brantingham, and G.~E. Tita.
\newblock {Measuring and modeling repeat and near-repeat burglary effects}.
\newblock {\em Journal of Quantitative Criminology}, 25(3):325--339, may 2009.

\bibitem{Short2008}
M.~B. Short, M.~R. D'Orsogna, V.~B. Pasour, G.~E. Tita, P.~J. Brantingham, A.L.
  Bertozzi, and L.~B. Chayes.
\newblock {A statistical model of criminal behavior}.
\newblock {\em Math. Models Methods Appl. Sci.}, 18(Suppl.):1249--1267, 2008.

\bibitem{Smith2012}
L.~M. Smith, A.~L. Bertozzi, P.~J. Brantingham, G.~E. Tita, and M.~Valasik.
\newblock {Adaptation of an animal territory model to street gang spatial
  patterns in Los Angeles}.
\newblock {\em DCDS}, 32(9):3223--3244, 2012.

\bibitem{stinner_win}
C.~Stinner and M.~Winkler.
\newblock {Global weak solutions in a chemotaxis system with large singular
  sensitivity}.
\newblock {\em Nonlinear Analysis: Real World Applications}, 12(6):3727--3740,
  2011.

\bibitem{tao_small}
Y.~Tao.
\newblock {Boundedness in a chemotaxis model with oxygen consumption by
  bacteria}.
\newblock {\em Journal of Mathematical Analysis and Applications},
  381(2):521--529, 2011.

\bibitem{taowin_consumption}
Y.~Tao and M.~Winkler.
\newblock {Eventual smoothness and stabilization of large-data solutions in a
  three-dimensional chemotaxis system with consumption of chemoattractant}.
\newblock {\em Journal of Differential Equations}, 252(3):2520--2543, 2012.

\bibitem{Tse2008}
W.~H. Tse and M.~J. Ward.
\newblock {Hotspot formation and dynamics for a continuum model of urban
  crime}.
\newblock {\em European J. Appl. Math.}, pages 583--624, 2016.

\bibitem{win_MANA}
M.~Winkler.
\newblock {Absence of collapse in a parabolic chemotaxis system with
  signal-dependent sensitivity}.
\newblock {\em Mathematische Nachrichten}, 283(11):1664--1673, 2010.

\bibitem{win_JDE}
M.~Winkler.
\newblock {Aggregation vs. global diffusive behavior in the higher-dimensional
  Keller-Segel model}.
\newblock {\em Journal of Differential Equations}, 248(12):2889--2905, 2010.

\bibitem{win_chemosing}
M.~Winkler.
\newblock {Global solutions in a fully parabolic chemotaxis system with
  singular sensitivity}.
\newblock {\em Mathematical Methods in the Applied Sciences}, 34(2):176--190,
  2011.

\bibitem{win_JMPA}
M.~Winkler.
\newblock Finite-time blow-up in the higher-dimensional parabolic-parabolic
  keller-segel system.
\newblock {\em J.~Math.~Pures Appl.}, 100:748--767, 2013.

\bibitem{win_ARMA}
M.~Winkler.
\newblock {Stabilization in a two-dimensional chemotaxis-Navier-Stokes system}.
\newblock {\em Archive for Rational Mechanics and Analysis}, 211(2):455--487,
  2014.

\bibitem{win_ks_nasto_exist}
M.~Winkler.
\newblock {A three-dimensional Keller-Segel-Navier-Stokes system with logistic
  source: Global weak solutions and asymptotic stabilization.}
\newblock {\em Preprint}, 2016.

\bibitem{Zipkin2014}
J.~Zipkin, M.~B. Short, and A.L. Bertozzi.
\newblock {Cops on the dots in a mathematical model of urban crime and police
  response}.
\newblock {\em DCDS-B}, 19(0):1479--1506, 2014.

\end{thebibliography}

\end{document}